\pdfoutput=1
\documentclass{article}
\usepackage[margin=25mm]{geometry}
\usepackage{amsmath}
\usepackage{amsfonts}
\usepackage{amssymb}
\usepackage{graphicx}
\usepackage[utf8]{inputenc}
\usepackage[T2A]{fontenc}
\usepackage[ukrainian,russian,english]{babel}

\usepackage{xcolor}
\usepackage[all]{xy}

\newtheorem{theorem}{Theorem}

\title{Three-color graph as the  1-skeleton of the 2-sphere triangulation}
\author{Oleg Akchurin, Svitlana Bilun and Alexandr Prishlyak}

\begin{document}
\maketitle

\begin{abstract}
The paper is devoted to finding the colorings of the edges of the 1-skeleton of triangulations of the 2-sphere in three colors so that for each face all three of its sides have different colors. First, by the method of adding one vertex inside the triangle or on its side, we enumerate all tiangulations with no more than 8 vertices. Next, one triangulation with 6 and 7 vertices, each with two different colors, was found. And finally, it is shown that other triangulations, which have less than 8 vertices, have one coloring each.
\end{abstract} \hspace{10pt}


\section*{Introduction}

A three-color graph is a convenient tool for classifying Morse functions and Morse-Smale vector fields on surfaces. For its construction, the division of the sphere into triangles is used. To do this, one regular trajectory is carried out in each cell that repeats the separatrices. Sepratrices entering saddles are colored in the first color (red), those coming out of saddles are colored green, and regular trajectories are colored blue. Thus, the resulting division into triangles has the property that the sides of each triangle are colored in different colors. A dual graph is called a tricolor graph of the Morse-Smale flow and the corresponding Morse function. Therefore, such graphs are in mutually unambiguous correspondence with the three-color coloring of the 1-skeleton of surface triangulation.

The aim of the paper is to find all possible colorings of triangulations of the 2-sphere with a small number of vertices.

It should be noted that, in addition to the tricolor graph, there are many other topological invariants for functions \cite{prishlyak2001conjugacy, hladysh2019simple, hladysh2017topology, prishlyak2003regular, prishlyak2002topological1, prishlyak2002morse, prishlyak2000conjugacy, prishlyak1999equivalence, prishlyak2007classification, lychak2009morse, prish2002Morse, prish2015top, prish1998sopr, bilun2002closed, Reeb1946, Kronrod1950, Sharko1993} and flows \cite{prislyak2017morse, prishlyak2019optimal, prishlyak2003topological, prishlyak2003sum,  prishlyak2005complete, prishlyak2020three, kkp2013, prish1998vek, prish2002vek, prish2001top, prishlyak2002morse1, prishlyak2002topological, Peixoto1959, Prishlyak2017, Prishlyak2022, Prishlyak2021, Hatamian2020, Prishlyak2020, Kybalko2018, Prishlyak2002, Prishlyak2019, Prishlyak2007, Peixoto1973, Palis1982, Palis1970, Fleitas1975, Palis1968, BPP2022, Giryk1996, Smale1961, Oshemkov1998, Bolsinov2004, Kadubovskyj2005, Poltavec1995, Morin1978, Andronov1937, Leontovich1955, Smale1960} on surfaces and 3-manifolds.
We also consider the triangulation (triangulation) of the surface itself as a graph embedded in the surface (its 1-skeleton). You can get acquainted with the main invariants of graphs and their embeddings in surfaces in \cite{prishlyak1997graphs, Harary69, pontr86, tatt88, HW68, GT87}.

Structure of this publication. In section 1 we describe all triangulations with no more than 6 vertices, in section 2 -- with 7 vertices, in section 3 -- with 8 vertices. In section 4, we find all possible colorings of the graphs described in sections 1 and 2.

To find new graphs, we add 1 vertex to already found graphs. At the same time, if the vertex is added to the middle of the triangle, three more edges are added along with it. If a vertex is added to an edge, it splits that edge into two and adds two more edges in addition. In this way, you can compile a list of all graphs containing less than 12 vertices. This follows from Euler's formula and the fact that each face contains three edges. In the resulting list, the graphs are checked for isomorphism. By ${\mu}( n )$ we denote the number of non-isomorphic graphs with $n$ vertices.

To construct colorings, we arbitrarily color one of the triangles (the outer one) and search for all possible continuations of this coloring.

\newpage

\section{Triangulations with 4, 5 and 6 vertexes}

Obviously, the minimum quantity { ie} n = 4, for any one ( ${\mu}$( n ) = 1) in the form of a split is a tetrahedron
(Fig. 1). The viable graph is meaningfully G \textsubscript{4 }. Respectfully, if instead of the task of splitting on
tricots consider the problem of triangulation, then the minimum number of vertices required for triangulation of the
sphere, obviously, is more than 3 (vertices are split along the equator, upper and lower pins are cut by
triangulations). Please note that, no matter how important it is to show, for the triangulation of a closed surface in
R \textsuperscript{3 }, a sphere, one point is sufficient (to bring it up, it is necessary to look at the process of
gluing a fairly large surface in R \textsuperscript{3 }from a polyhedron).

The task of splitting the sphere into tricks for n {\textgreater} 4 will be overcome by the offensive rank. Let us give
you a set of differences (to the point of isomorphism) of the division of the sphere for each n . Let's take a breakup
G \textsubscript{nk }. At first glance, it's possible to go from one to the other before breaking G \textsubscript{n
+1} two and only two ways are possible (see. fig. 1):

\begin{enumerate}
\item choosing the n +1-th vertex in the middle of the face (triangle) to the graph (decomposition ) G
\textsubscript{nk}
\item choosing the n +1 -th vertex on the next edge of the graph G \textsubscript{nk} opposing vertexes of similar adjacent triangles
\end{enumerate}
However, in this way, new vertices of step n =3 (for slope 1) and step n =4 (for slope 2) are added. Obviously, the
procedure (1) and (2) is by no means obsessed with the icosahedron, all the vertices of which are able to make steps 5.
The reason for this is the need for a search for less trivial options for the transition from the group G
\textsubscript{nk }to the graph G \textsubscript{n +1 }. Can you show that for whom it is sufficient to supplement the
procedures (1) and (2) with an offensive procedure (see Fig. 2):

3) Choose in the graph G \textsubscript{nk }a real n '-gon that does not have an internal vertex.

Anuluvati 2 inner ribs of the p ' -gon. I will step on the n + 1st peak in the middle of the p '-gon. 3 ' one newly formed vertex from all vertices n ' -gon.

Bringing the sufficiency of listing procedures (1)-(3) to encourage any kind of splitting of the sphere into triangles is
evident from the theory of power graphs about those that have a planar graph skin (and these graphs are equivalents of
splitting, obviously, planar) have the vertex of the degree not more than 5 \cite{Harary69, tatt88}.

\begin{figure}[ht]  
\center{\includegraphics[height=6.5cm]{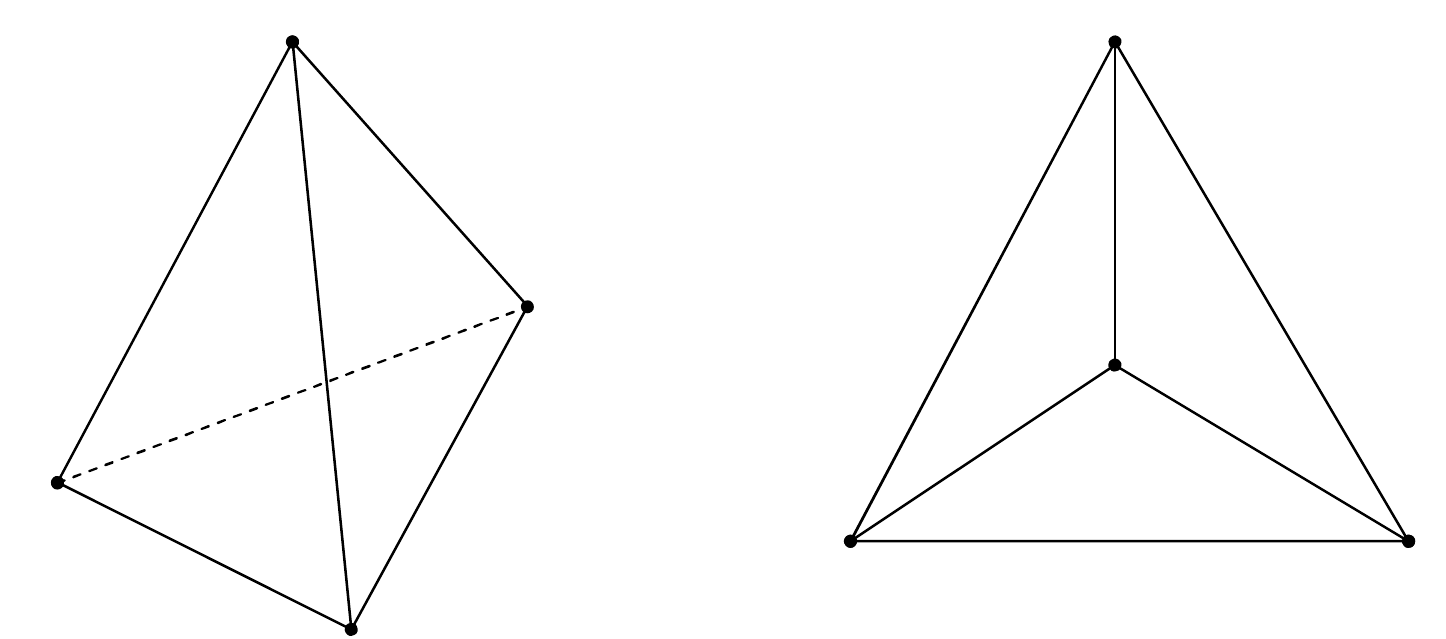}}
\caption{Breaking the sphere into triangles with n =4, ${\mu}$( n )=1.} 
\label{p1c} 
\end{figure}

It was conceivable that in such a way all ${\mu}$( n +1) edges on triangles of the sphere with an additional n +1
vertices will be won. However, acts from them, naturally, can appear isomorphic. In this order, the task is reduced to
re-verification of the isomorphism of obsessions as a result of the described procedure G \textsubscript{n +1 }{}-
graphs.

\begin{figure}[ht] 
\center{\includegraphics[height=6.5cm]{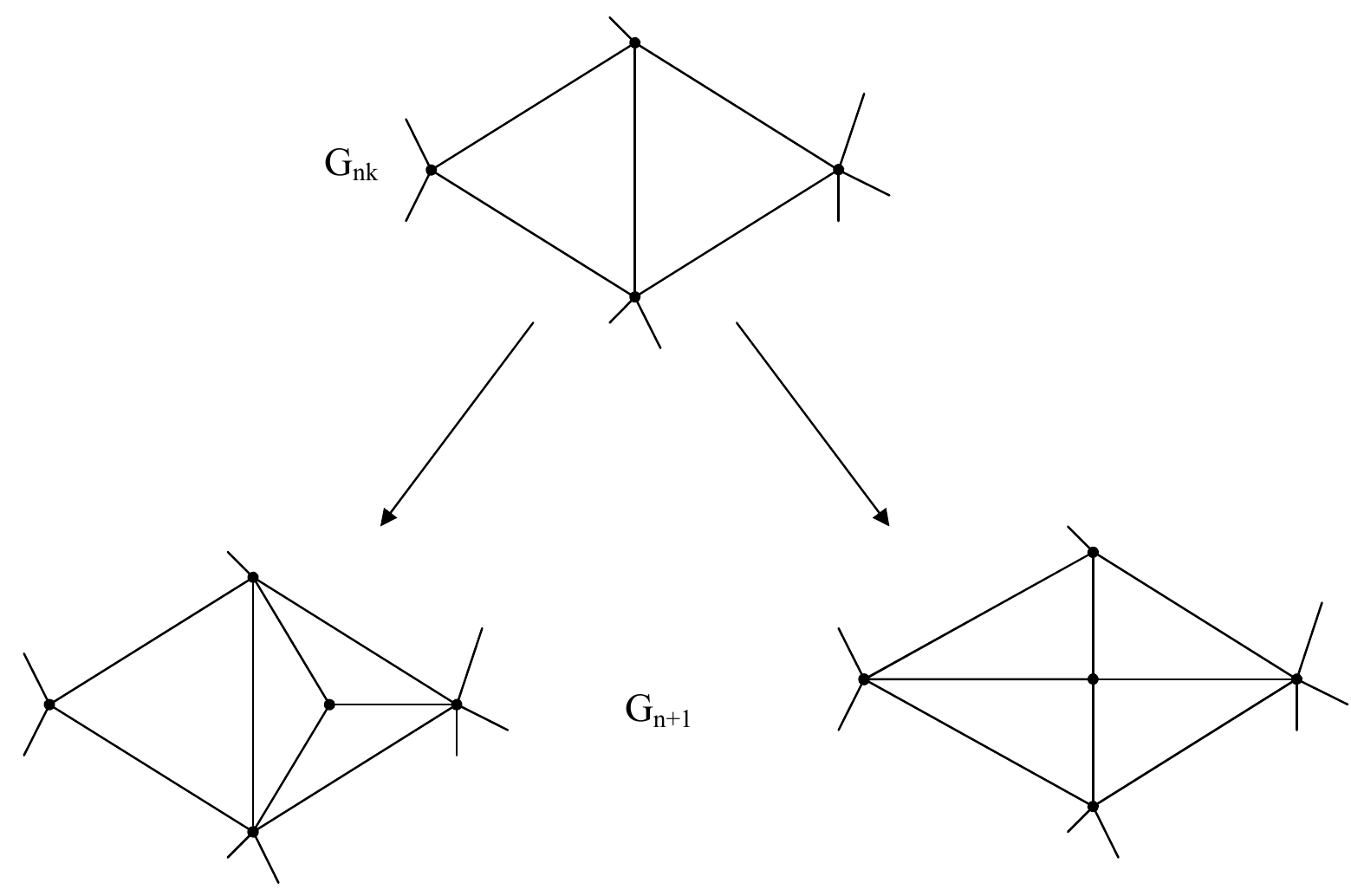}}

\caption{$ $} 
\label{p2} 
\end{figure}

Just one more reason to respect them, how to go without a middle until you break the spheres into tricks for specific n. In general, the procedure for the transition from G \textsubscript{nk } to G \textsubscript{n +1 } may seem rather
laborious. It's possible to make it easier, taking care of your respect, as if it 's possible , the manifestation of
symmetry.
Possible options (1 - left-handed) and (2 - right-handed) put the graph G \textsubscript{n +1 } behind
the graph G \textsubscript{nk } (Fig. \ref{p2})

\begin{figure}[ht]  \center{\includegraphics[height=6.5cm]{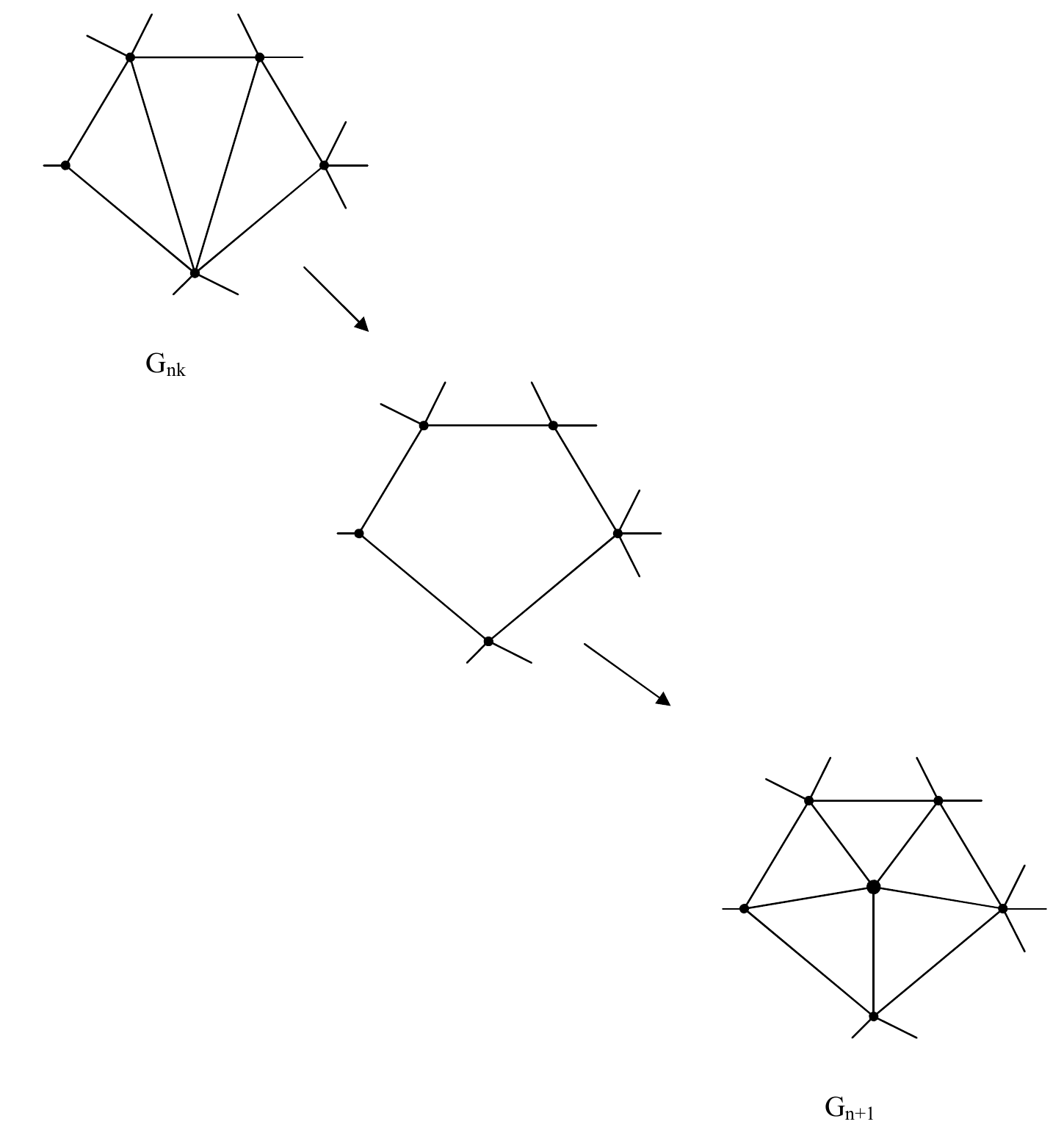}}
\caption{$ $} 
\label{p3} 
\end{figure}

Passing to the tetrahedron G\textsubscript{4 } to splitting with additional 5 points can be taken to respect that all
faces of the tetrahedron are absolutely symmetrical; For this, there is no need to carry out a fight for all sides, but
it's enough to settle down with only one; similar for ribs. In this order, there are two possible divisions of the
sphere with additional 5 points (Fig. \ref{p3}); it does not matter to show that the stench is isomorphic. In this order,
${\mu}$(5)=1.

In option (3) put the graph G \textsubscript{n +1 } behind the graph G \textsubscript{nk }. (Fig.\ref{p3})

\begin{figure}[ht]  \center{\includegraphics[height=6.5cm]{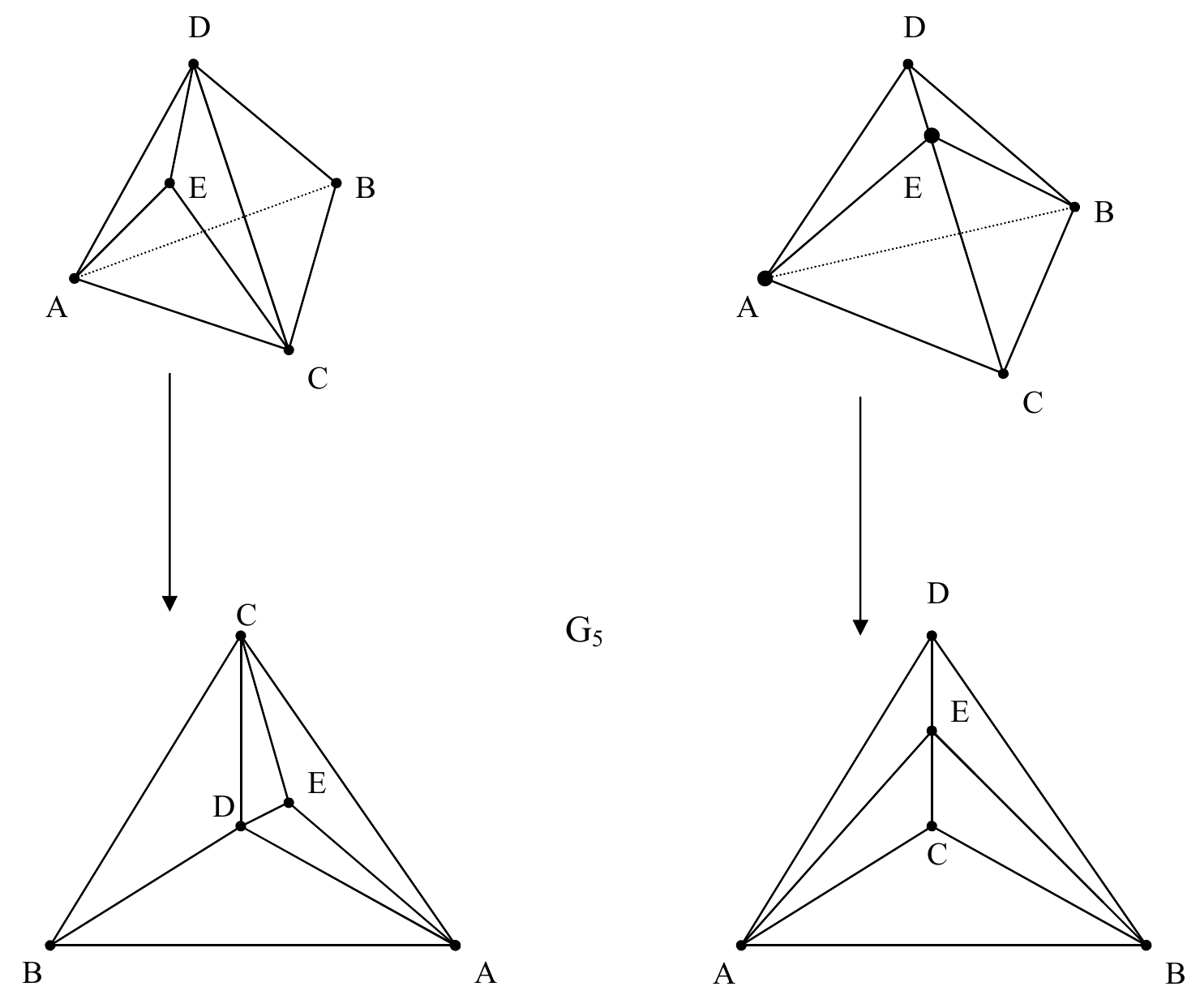}}
\caption{Graph $G_5$} 
\label{p4} 
\end{figure}

Variants of cutting spheres on triangles with additional 5 points. The graphs are left-handed and
right-handed isomorphic: A \textsubscript{l }\~{} A \textsubscript{p }, C \textsubscript{l }\~{} B \textsubscript{p },
E \textsubscript{l }\~{} C \textsubscript{p }, D \textsubscript{l }\~{} E \textsubscript{p }, B \textsubscript{l }\~{}
E \textsubscript{p } (Fig. \ref{p4})

\begin{figure}[ht]  \center{\includegraphics[height=8.5cm]{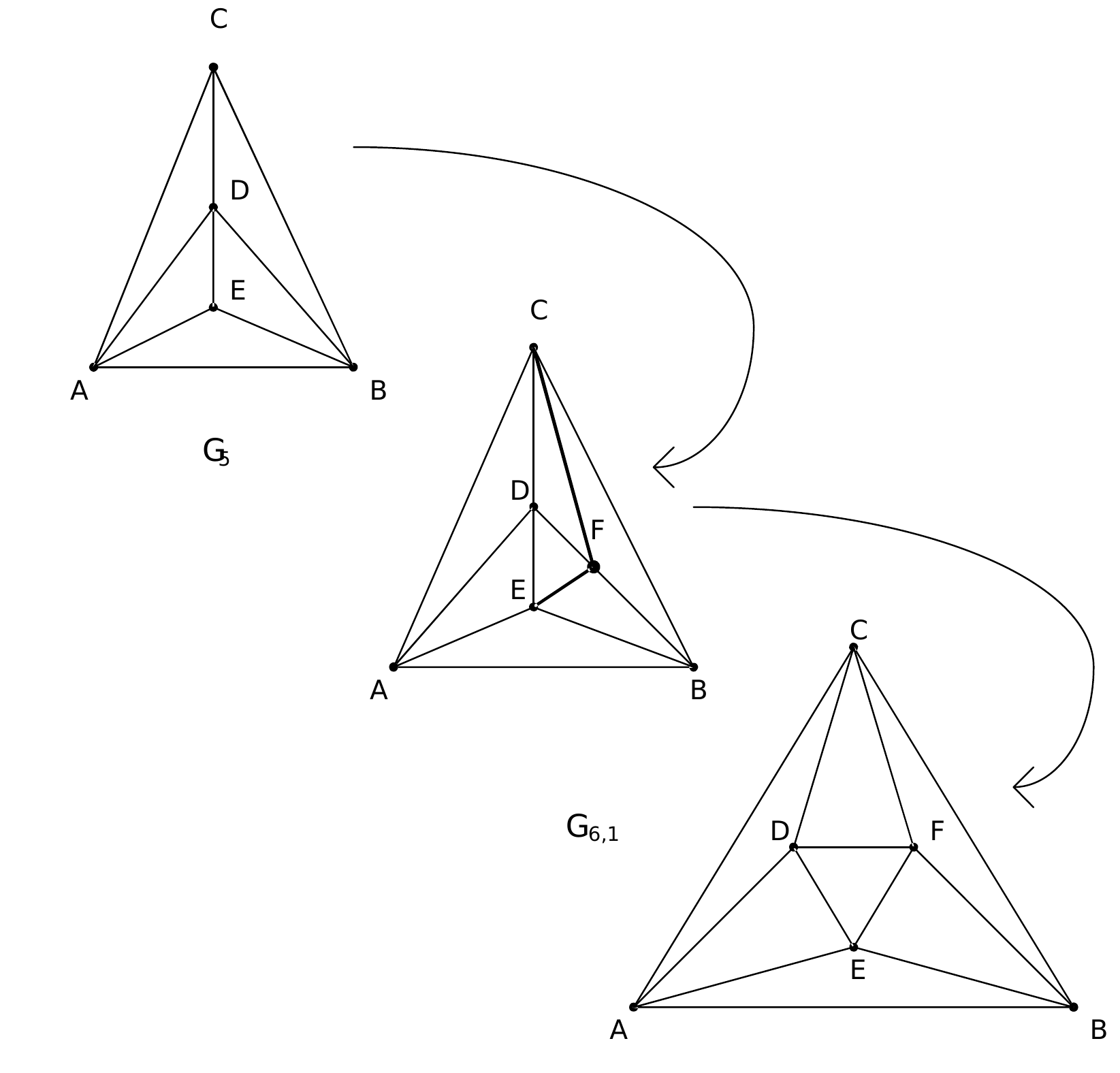}}\caption{Graph $G_{6,1}$} \label{p5} \end{figure}

Let's move on to looking at the possessed graph G \textsubscript{5 }. It is not important to remember that geometrically
you form a figure, as it is two ``glued'' tetrahedrons along the face. Obviously, all the faces of G
\textsubscript{5} are symmetrical. If there are edges, then, for example, the edge AB (Fig. 5, left hand) is not
symmetrical to the edge AC, so the steps of the vertices of the edge AB are equal to 3 and 4, and the edges AC - 4 and
4. all the ribs to the graph G \textsubscript{5 are selected by options }. In such a way, in order to go to the
division of the sphere to tricots with additional 6 points, it is necessary to make 3 divisions of the graph G
\textsubscript{5 }, wrapping the 6th vertex on ``4-4''-ribs (1) (div. Fig. \ref{p5}), on the ``3-4''-ribs (2) and in the
middle, whether it be a face G \textsubscript{5 }(3) (div. Fig. 9) and check the possession of the graph for
isomorphism.

Breaking the sphere with 6 points: choose the vertex F on the ``3-4''-edges EB to the graph G
\textsubscript{5 }(2) and in the middle of the face AEB (3). The isomorphism of the possessed graphs is obvious (Fig. \ref{p5})

It is not important to check that options (2) and (3) give isomorphic graphs G 6,2 , and \textsubscript{option }(1) -
graph G 6,1 , which is not \textsubscript{isomorphic }to graph G \textsubscript{6,2 }. The non-isomorphism is evident
even if all the vertices of the graph G \textsubscript{6,1 }(octahedron) can have steps 4, just as the graph can have
vertices of steps 3, 4 and 5. Thus, $\mu $(6)=2.

\section{Triangulation with 7 vertexes}

Passing to the division of the sphere into triangles with the help of 7 points, next respect the step. First of all, we
can see the difference between 2 graphs - G \textsubscript{6.1 }and G \textsubscript{6.2 }. For skin problems, there
are traces of asymmetric elements and changes in the appearance of their differences. Let's look at the back of the
graph G \textsubscript{6.1 }. Osk{ ii}lki v{ ii}n v{ ii}dpovid{ a}{ ie} regular polyhedron - octahedron,
then in this case the situation is similar to { ii}{ ie}{ yi}, yak mimali for vipadku tetrahedron G
\textsubscript{4 }. And for us, it's enough for us to choose the 7th vertex on the  face (1) and on the
edge (2) graph G \textsubscript{6,1 }, even if the edges and vertices of the octahedron are symmetrical.

Looking at the graph G \textsubscript{6,2 }, remember (Fig. 9) that you can marvel at something like a tetrahedron, one
face of something ( CD ) is divided by two points ( E and F ) into 3 ribs, from the broken adjacent faces. From such a
look at G \textsubscript{6.2 }, we understand that there are 2 planes of symmetry - C DM and ABN ( M {}- the middle of
the edge AB, N {}- the middle of the edge CD ). Therefore, when passing from G \textsubscript{6,2 }to 7-folded spheres
on triangles, it is enough for us to look at the options for choosing 7-points on quiet elements, as they are found in
one of the quadrants, which are named symmetry planes.

In this rank, it is enough for us to look at the next options for choosing the 7th point on graph G \textsubscript{6.2
}(Fig. \ref{p6}):

\begin{figure}[ht]  \center{\includegraphics[height=14.5cm]{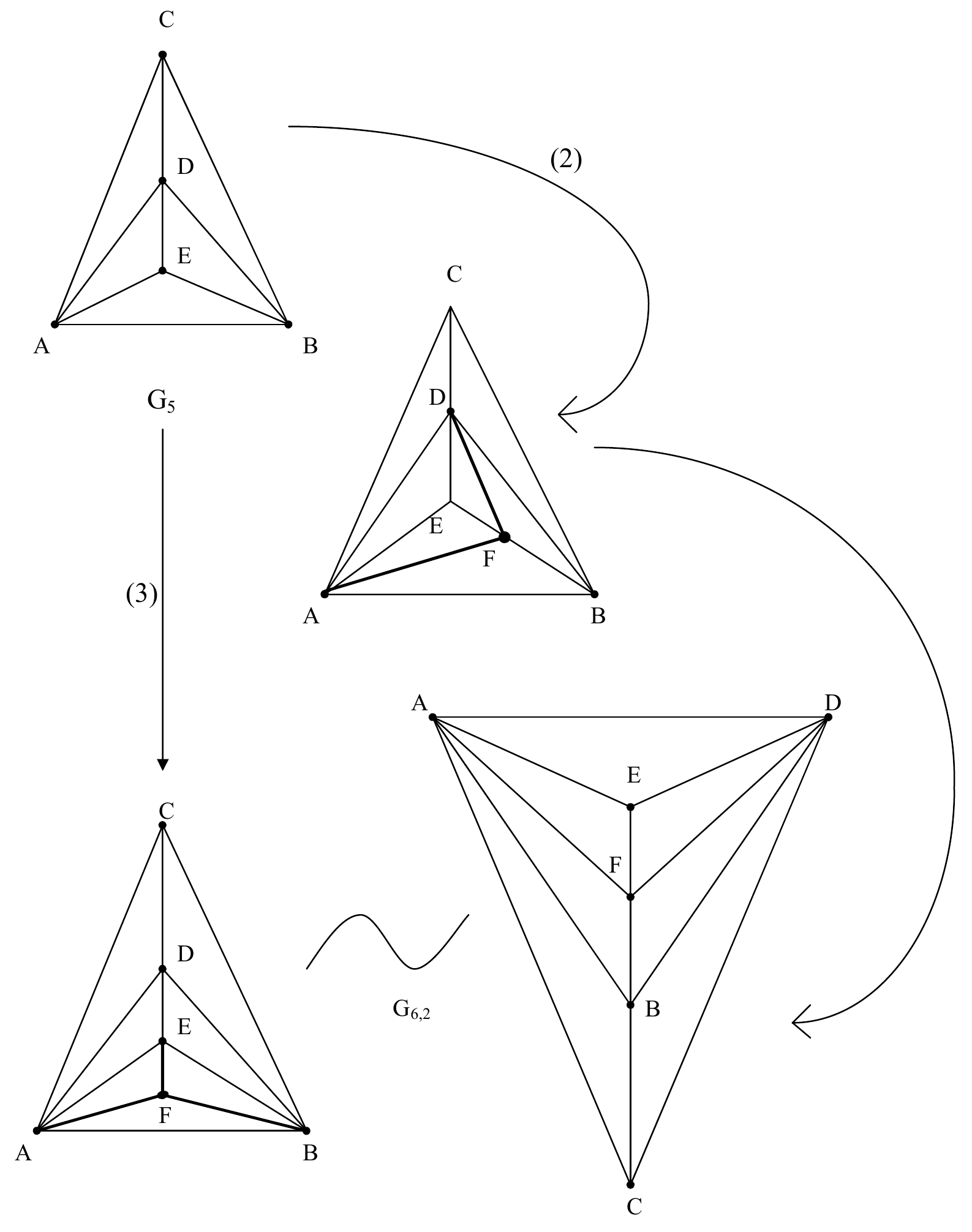}}
\caption{Graph $G_{6,2}$} 
\label{p6} 
\end{figure}

(3) -- on ``3-5''-ribs ({ S} A );

(4) -- on ``3-4''-ribs ({ S} F );

(5) -- on ``4-4''-ribs ( EF );

(6) -- on ``4-5''-ribs ( FA );

(7) -- on ``5-5''-ribs ( AB );

(8) - on `` 3-5-5 '' -faces (CAB);

(9) -- on ``3-4-5''-face ({ S} FA );

(10) -- on ``4-4-5''-face ( EFA ).

Let's take a closer look at the leather from the options (1) - (10).

\begin{figure}[t]  \center{\includegraphics[height=6.5cm]{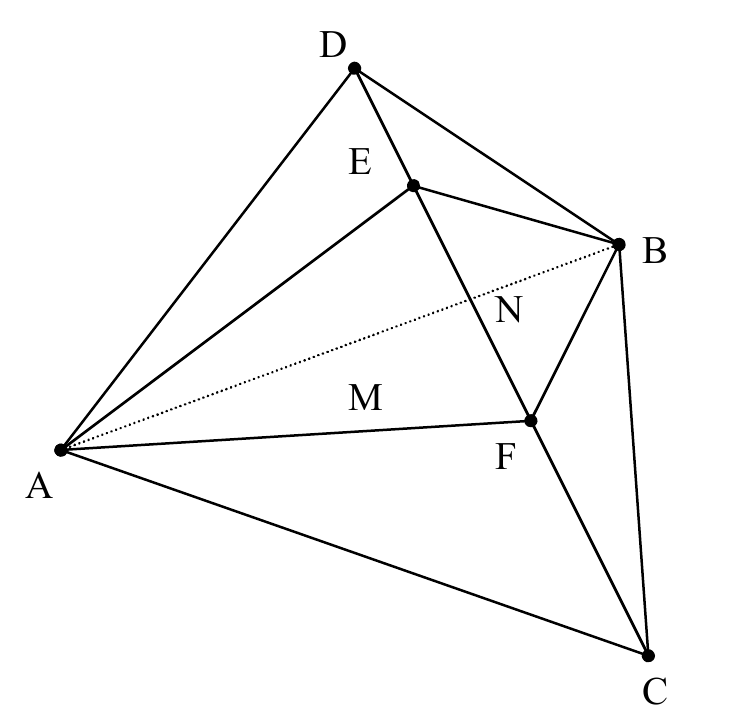}}
\caption{Graph G\textsubscript{6.2 } and other symmetry planes: DCM that ABN} 
\label{p8} 
\end{figure}
\hfill

\begin{figure}[ht]  
\center{\includegraphics[height=6.5cm]{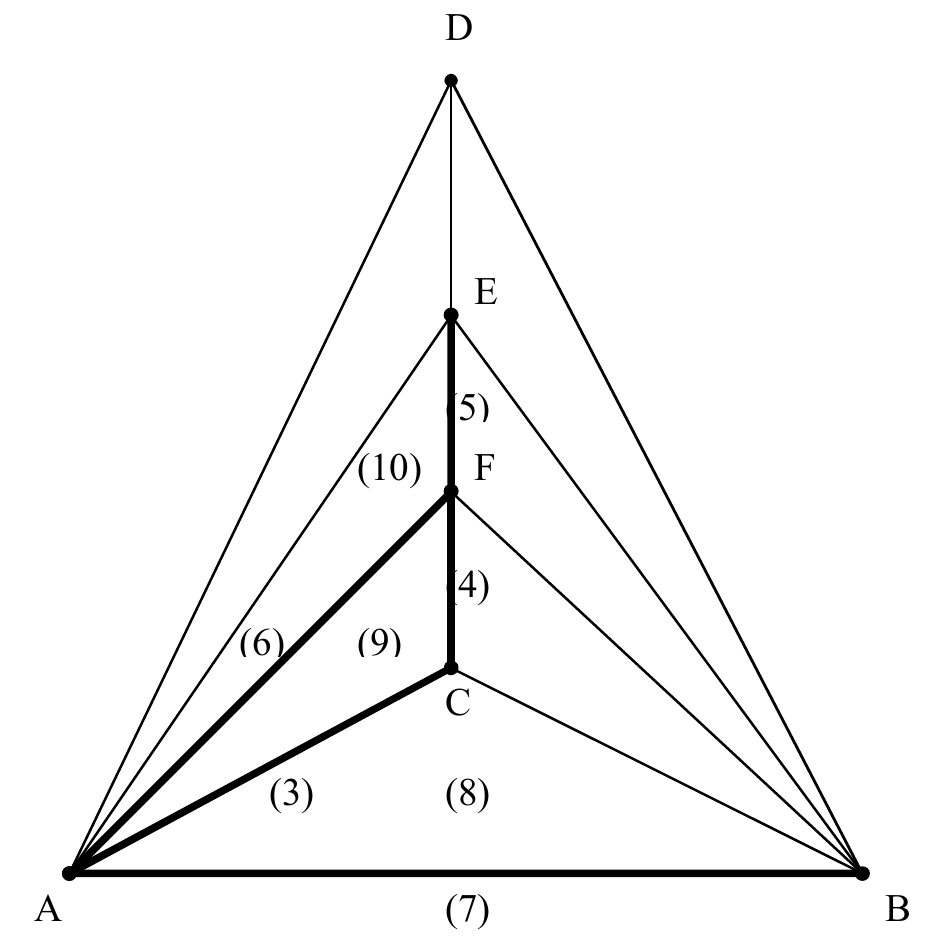}}
\caption{Graph $G_{6}$ vertexes. Elements were seen ((3) - (10)), on which sled to collect the 7th point to induce the splitting of the sphere into triangles with the help of 7 vertices.} 
\label{p8a} 
\end{figure}

\newpage

\begin{figure}[ht]  
\center{\includegraphics[height=16.5cm]{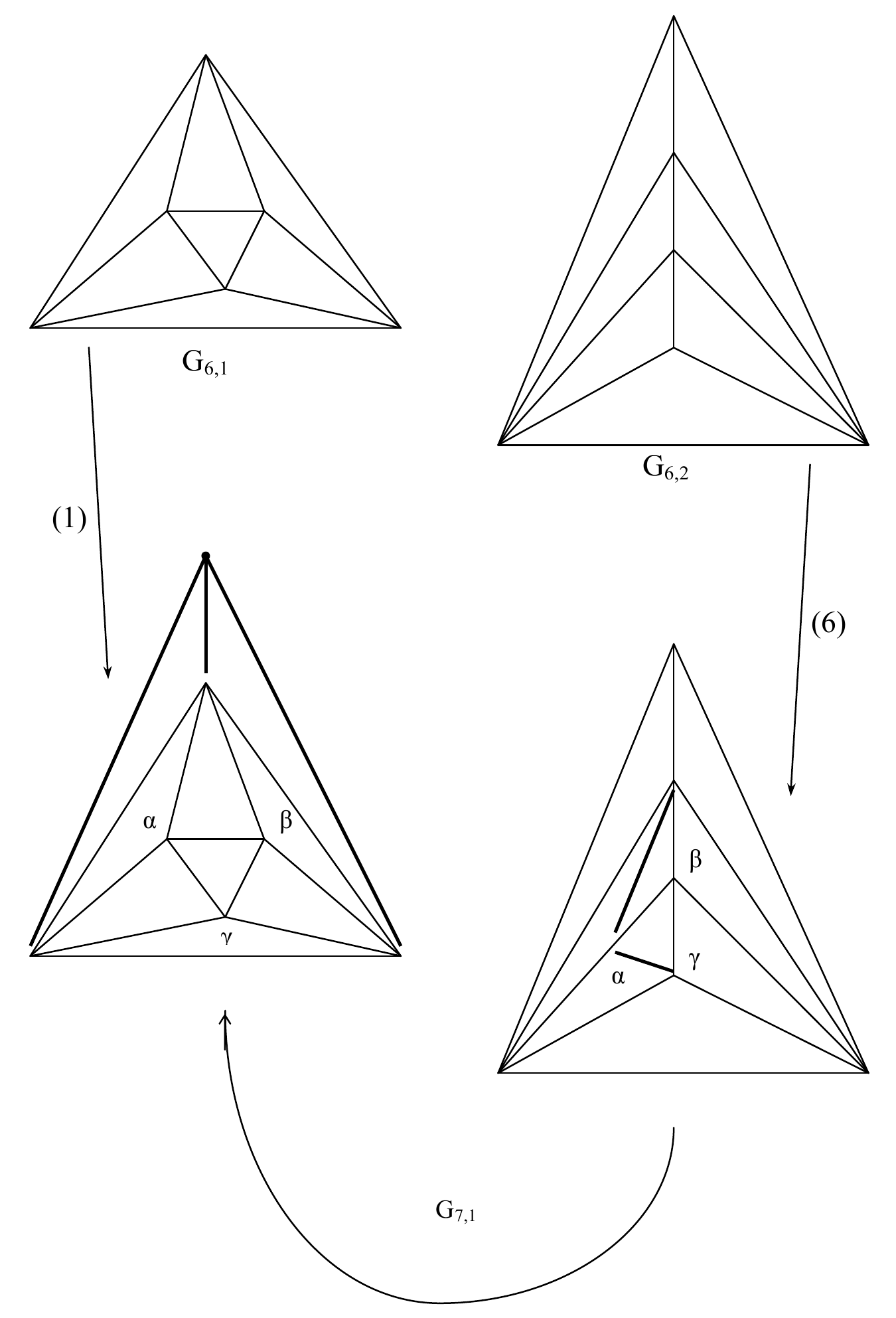}}
\caption{Graph $G_{7,1} $ edges. Variants (1) and (6) splitting the sphere into triangles with the help of 7 points. Obvious isomorphism of obsessive graphs.} 
\label{p9a} 
\end{figure}

\newpage

\begin{figure}[ht]  \center{\includegraphics[height=16.5cm]{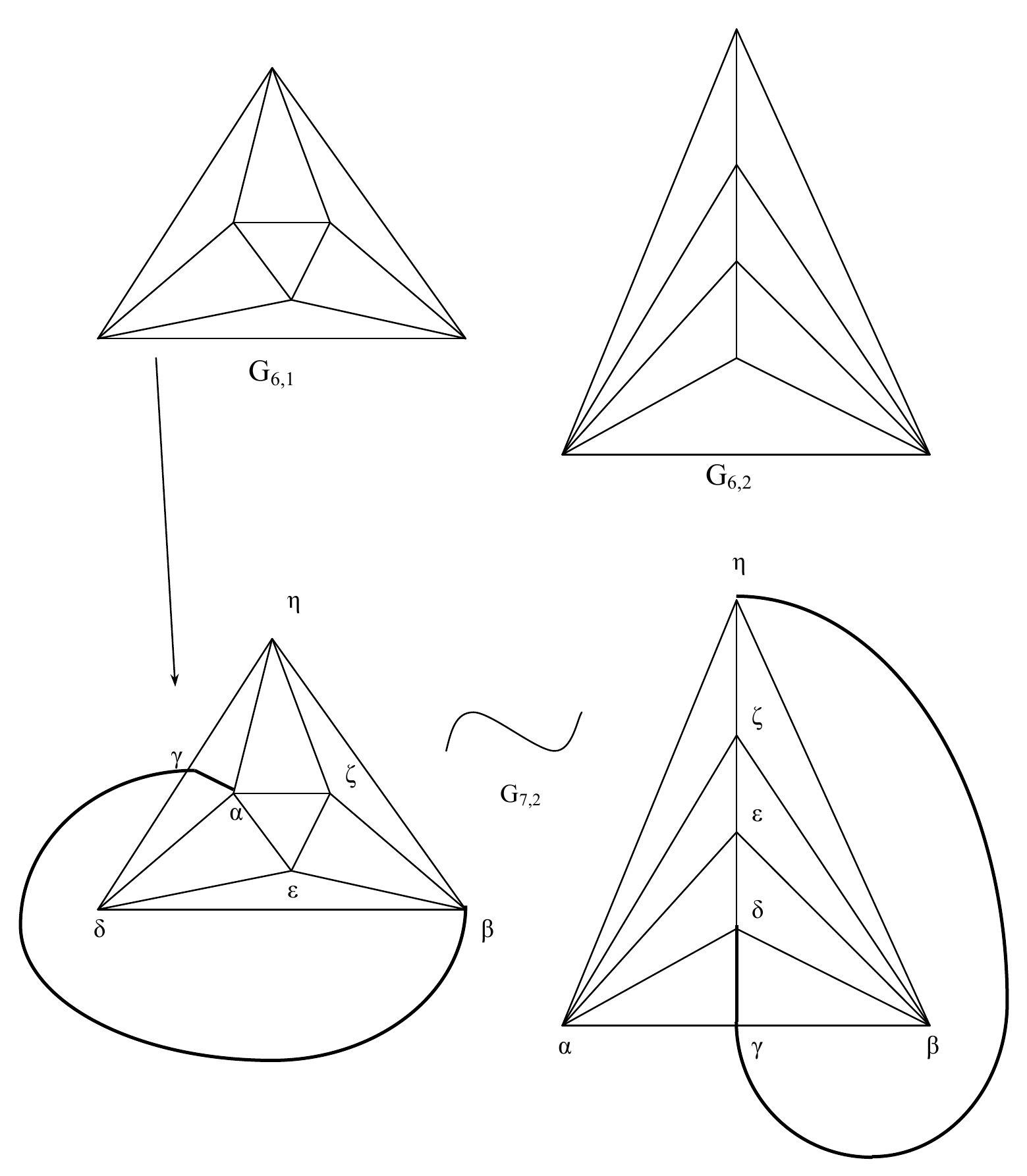}}
\caption{Graph $G_{7,2}$. Variant (2) and (7) splitting the sphere into triangles with the help of 7 points. The possession of
graphics isomorphic. Isomorphism is given to the same values of similar vertices. } 
\label{p10} 
\end{figure}

\newpage

\begin{figure}[ht]  \center{\includegraphics[height=15.5cm]{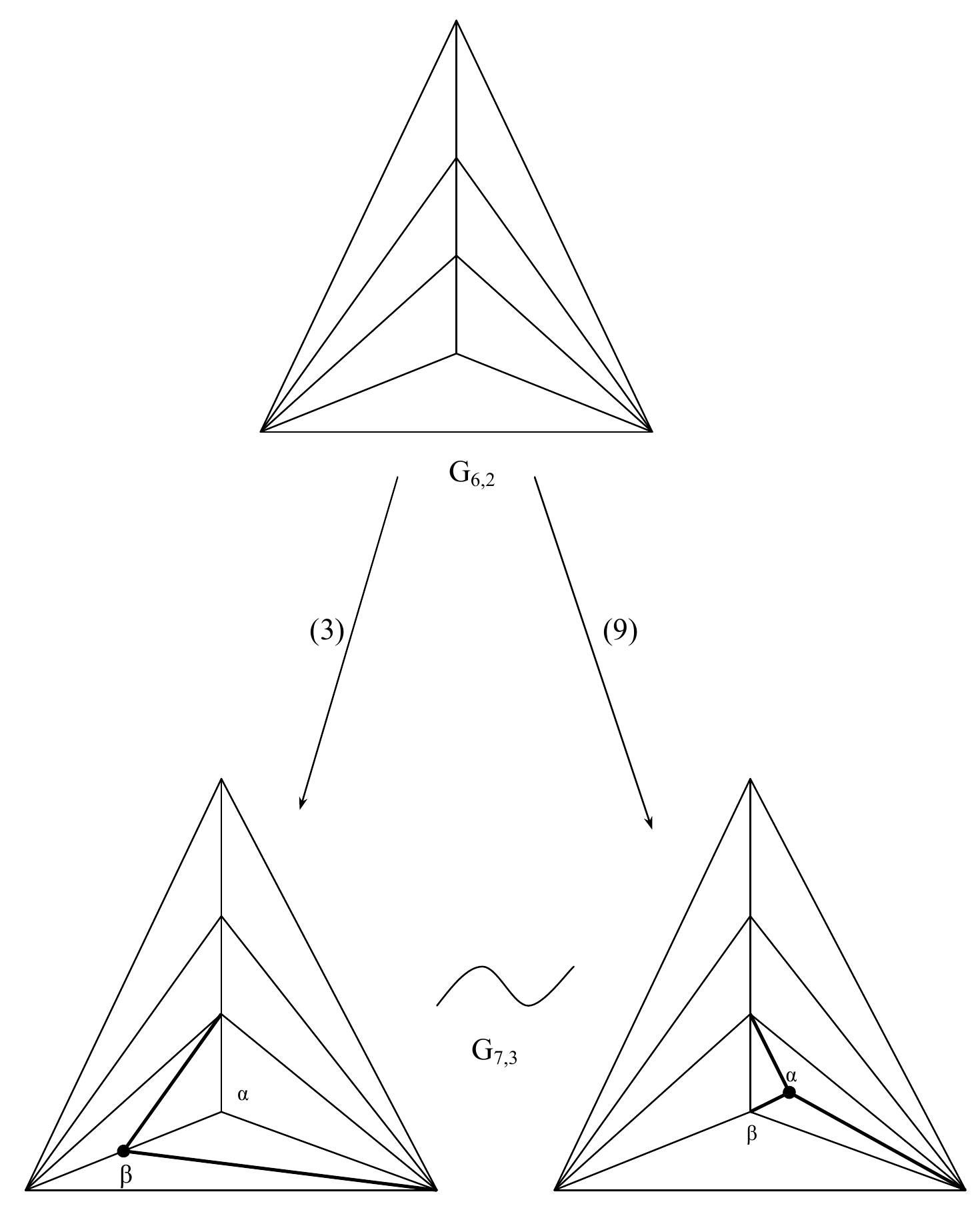}}
\caption{Graph $G_{7,3}$. Options ( 3 ) and ( 9 ) splitting the sphere into triangles with the help of 7 points. Obvious
isomorphism of obsessive graphs.} 
\label{p11} 
\end{figure}

\newpage

\begin{figure}[ht]  \center{\includegraphics[height=15.5cm]{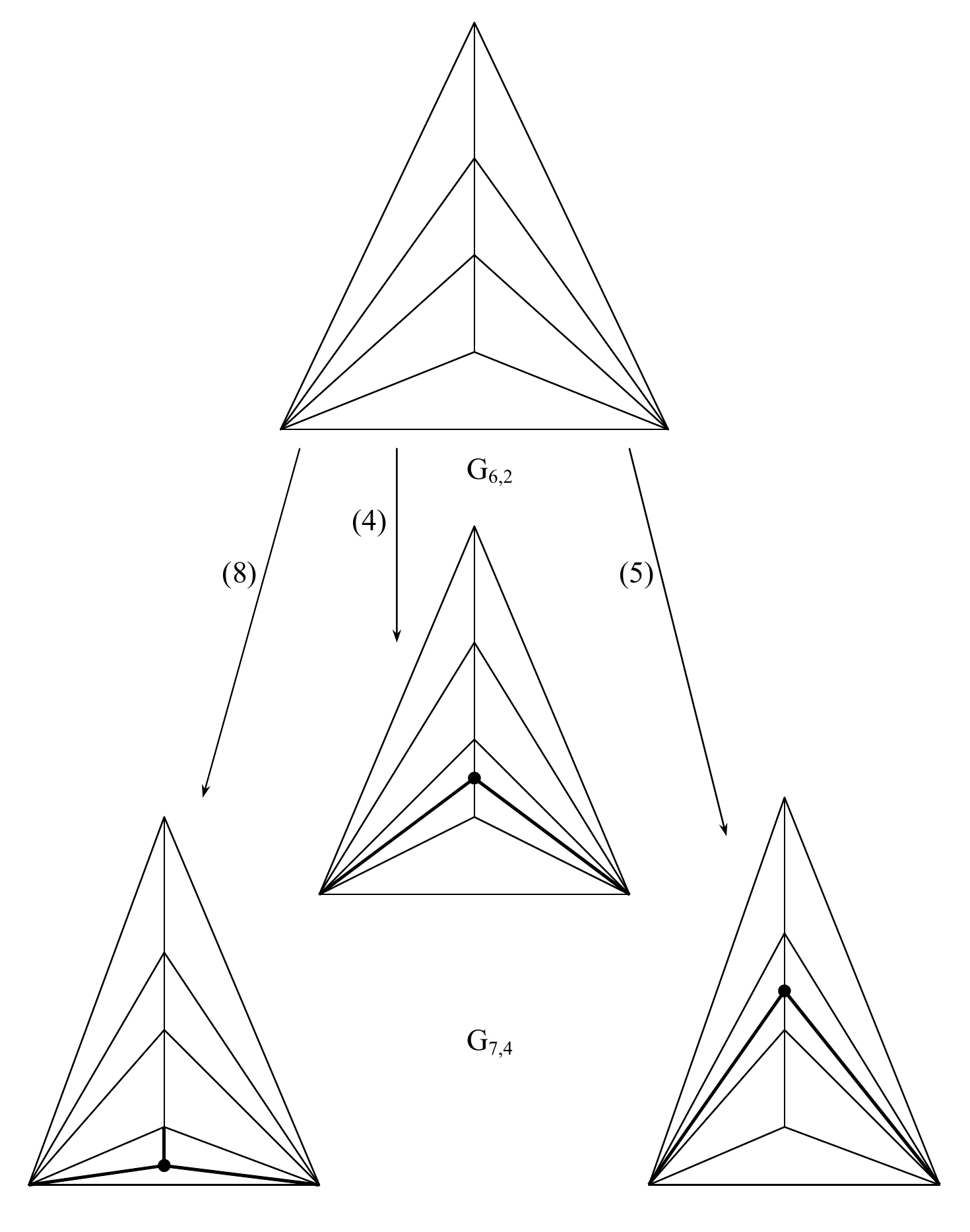}}
\caption{Graph $G_{7,4}$. Options ( 4 ) , (5) and ( 8 ) splitting the sphere into triangles with the help of 7 points. Obvious
isomorphism of obsessive graphs.} 
\label{p} 
\end{figure}

\newpage

\begin{figure}[ht]  \center{\includegraphics[height=11.5cm]{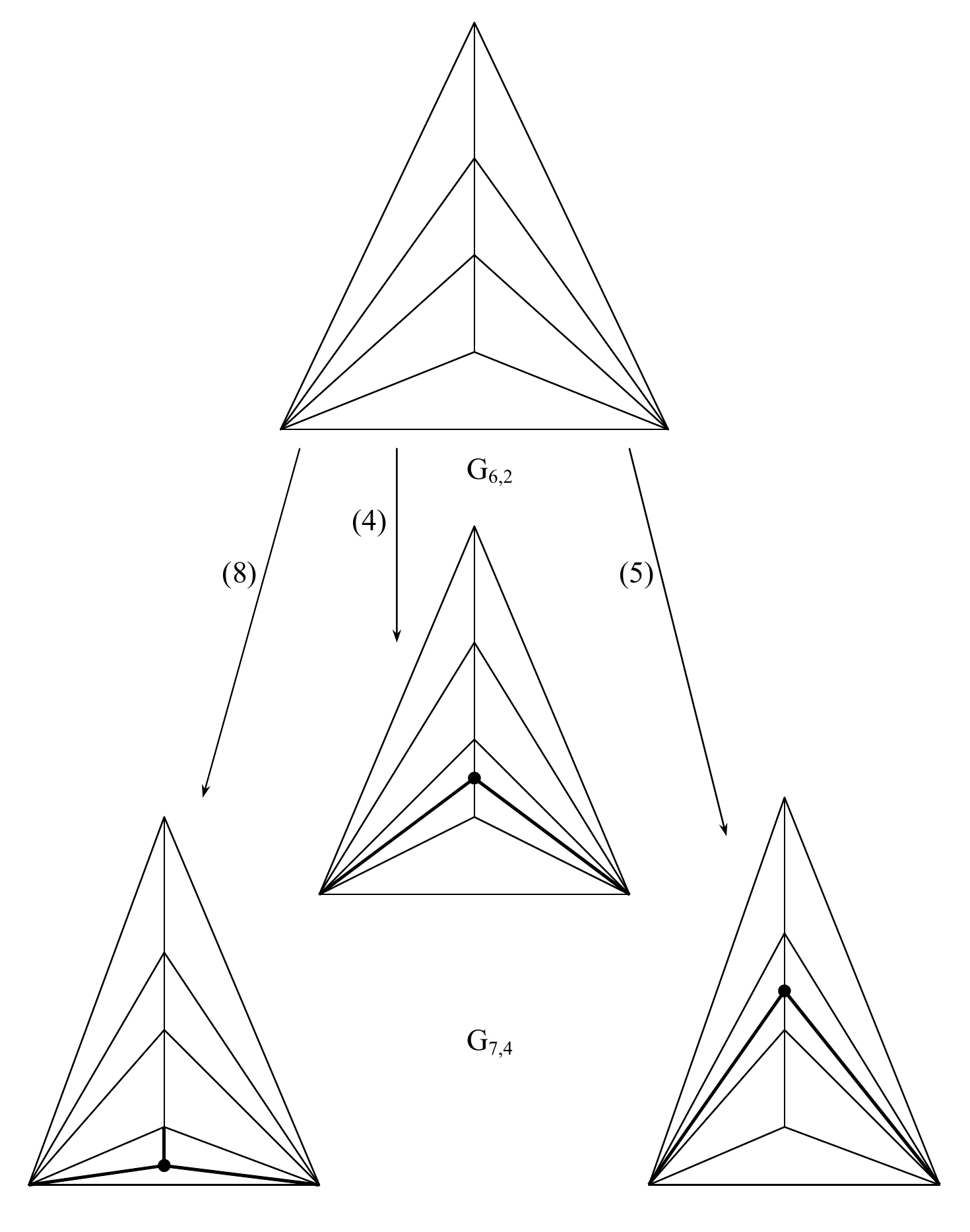}}
\caption{Graph $G_{7,5}$. Option (10) splitting the sphere into triangles with the help of 7 points.} 
\label{p13} 
\end{figure}

Are the graphs G \textsubscript{7.1 },{\dots}, G \textsubscript{7.5 }isomorphic? Table \ref{t1} gives a negative effect on the
power supply, in which a change of steps of the tops of these graphs is introduced.

\begin{table}[ht]
	\centering
		\begin{tabular} {|c|c|c|c|c|c|c|c|}
		
\hline
 G \textsubscript{7.1} &
 5 &
 5 &
 5 &
4 &
4 &
4 &
3 \\
\hline
 G \textsubscript{7.2} & 
\centering 5 &
\centering 5 &
\centering 4 &
\centering 4 &
\centering 4 &
\centering 4 &
 4 \\ \hline
\centering G \textsubscript{7.3} &
\centering 6 &
\centering 5 &
\centering 5 &
\centering 4 &
\centering 4 &
\centering 3 &
 3 \\ \hline
\centering G \textsubscript{7.4} &
\centering 6 &
\centering 6 &
\centering 4 &
\centering 4 &
\centering 4 &
\centering 3 &
 3 \\ \hline
\centering G \textsubscript{7.5} &
\centering 6 &
\centering 5 &
\centering 5 &
\centering 5 &
\centering 3 &
\centering 3 &
 3\\ \hline
\end{tabular}
\caption{Degree of vertexes} 
\label{t1} 
\end{table}

\bigskip

In this order, $\mu $(7)=5.

\newpage

\section{Triangulation with 8 vertexes}

Let's move on to the challenge of breaking out of 8 vertices, taking, I'll do it again, as a basis for winning 5
breaking out of 7 vertices. Respectfully, here we can't get along without the procedure (3) (on the front crochet, it's
not important to cross over, it doesn't give a welcome break, it's exactly the same kind of 5 wake-ups).

\begin{figure}[ht]  \center{\includegraphics[height=6.5cm]{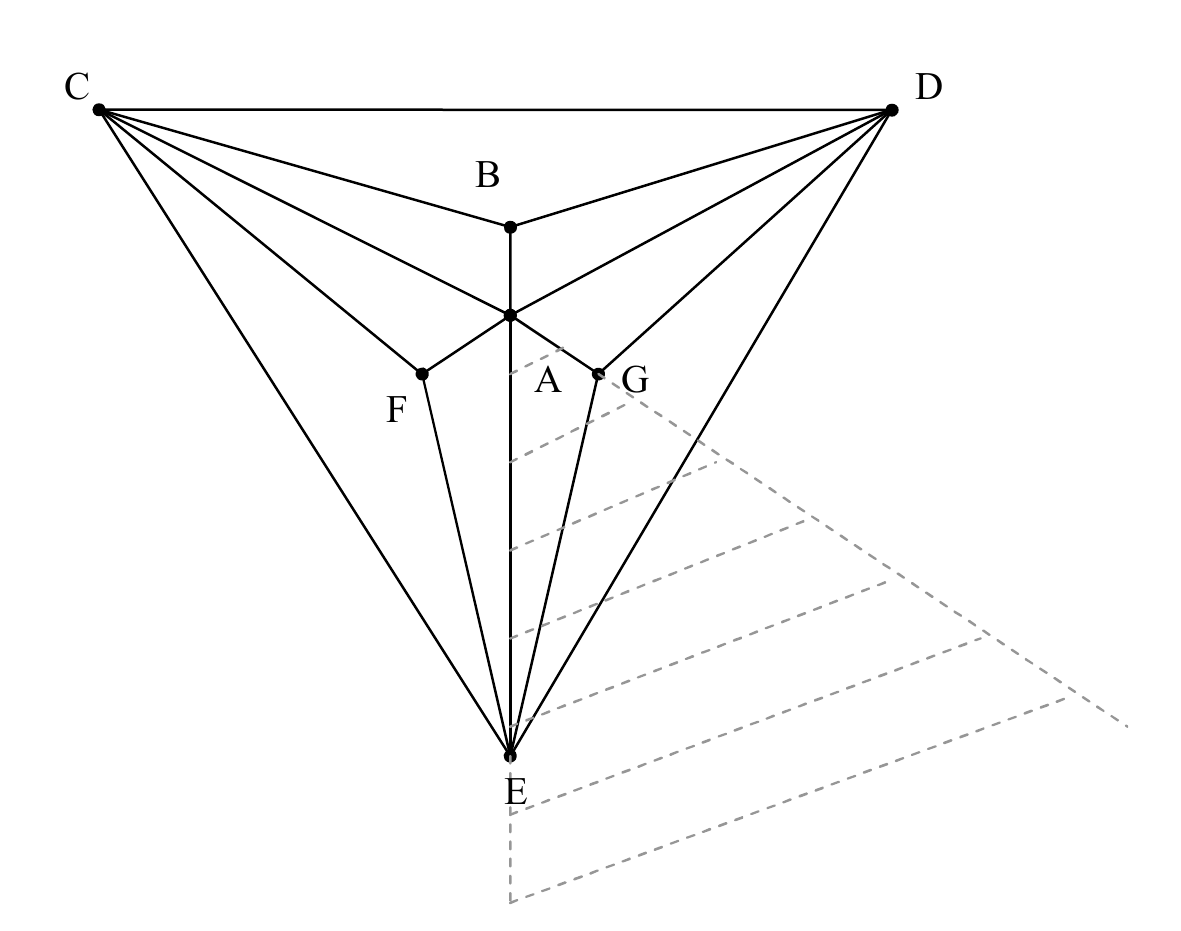}}\caption{Graph  G \textsubscript{7.5 }. The plot is shaded, as if looking at the distance of the battles, emerging from the world of symmetry.} 
\label{p14} 
\end{figure}

Let's take a look at graph G \textsubscript{7} in order, to the virtuous of their promptings and assignments by us of
numbering. Let's take a look at graph G \textsubscript{7.5 }. In Persian black, next to yoga in a more symmetrical look
(Fig. \ref{p14}). The isomorphism of the output graph is obvious -- it is sufficient to refer to the values of the isomorphic
vertices. In this order , it is obvious that in this time it is sufficient to look at only the breakdown of quiet
elements, which have reached the shaded area. Analyzing the number of elements, we can reach the conclusion, which should
look at the breakdown of the offensive and less offensive elements:

(1) AGE;

(2) EGD;

(3) CDE;

(4) ED;

(5) EG;

(6) EA ;

(7) AG ;

(8) EGAFC ;

(9) EGABC ,

where the two remaining battles mean the stagnation of reworking (3) until the last pentagon. Doing correspondent transformation, we obtain 19-24.

\newpage

\begin{figure}[ht]  \center{\includegraphics[height=7.5cm]{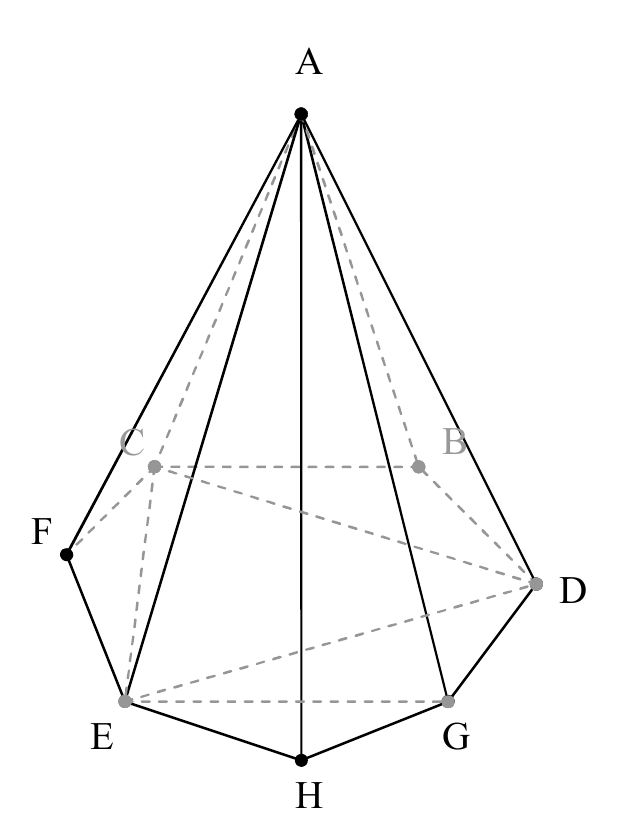}
\includegraphics[height=7.5cm]{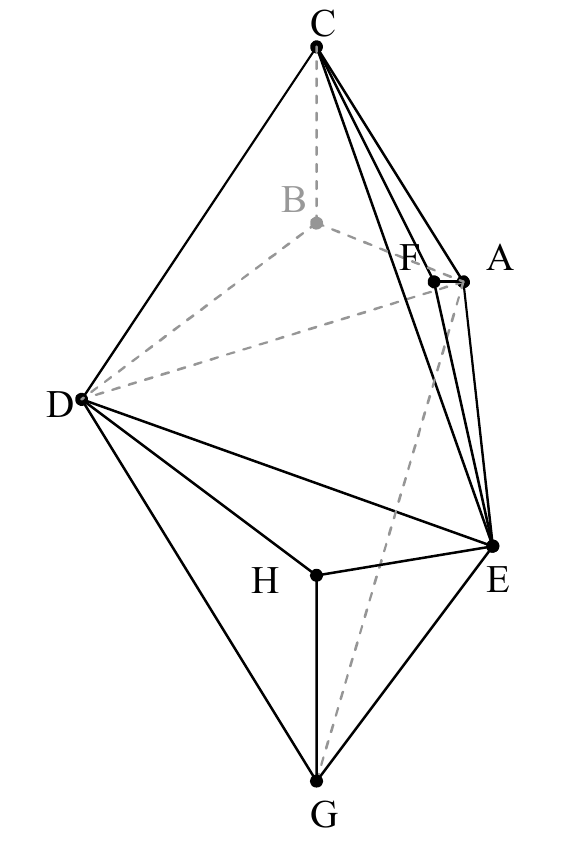}}\caption{Graphs $G_{8,1}$ and $G_{8,5}$} \label{p16} \end{figure}

\begin{figure}[ht] 
\center{\includegraphics[height=7.5cm]{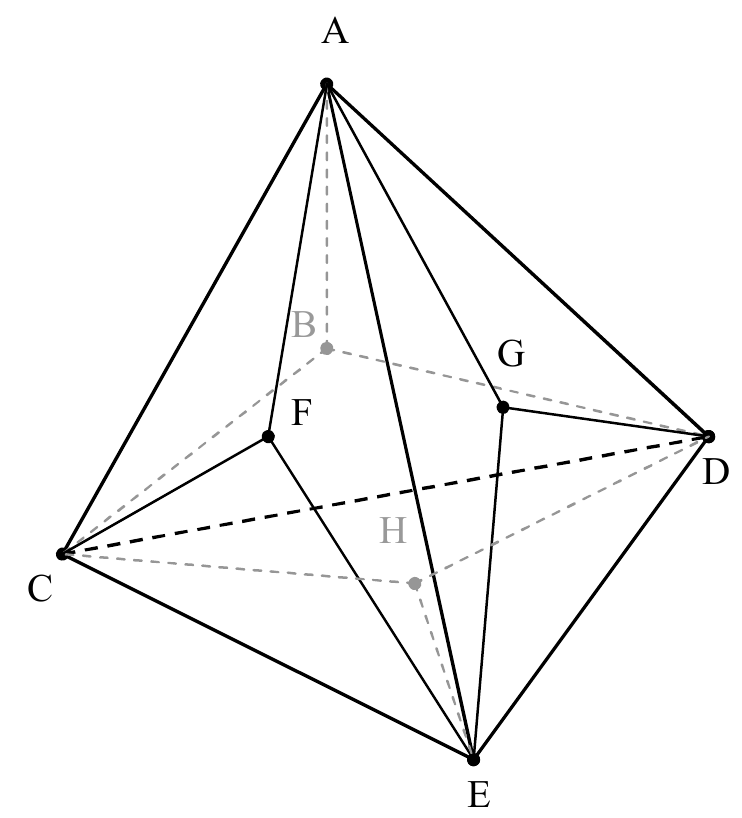}
\includegraphics[height=7.5cm]{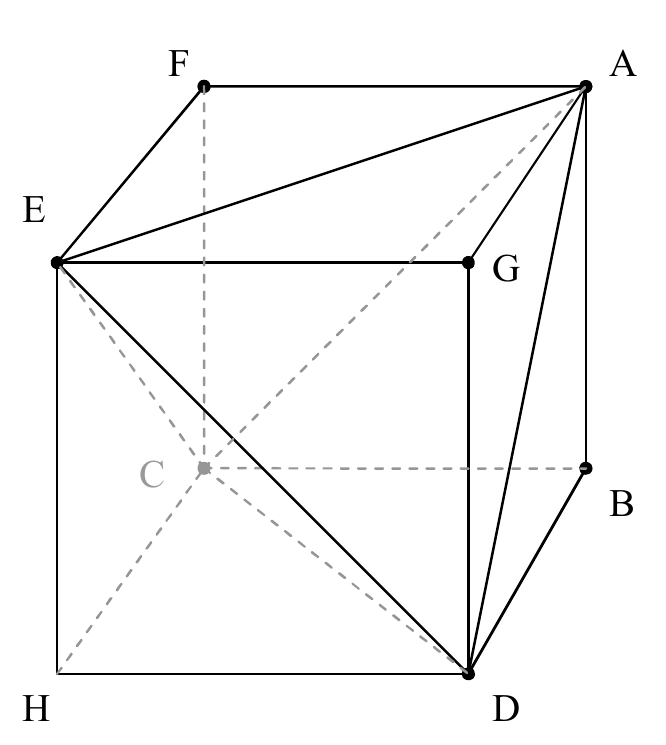}}\caption{Graphs $G_{8,2}$ and $G_{8,7}$} \label{p18} \end{figure}


\hfill
\newpage

\begin{figure}[ht]  \center{\includegraphics[height=7.5cm]{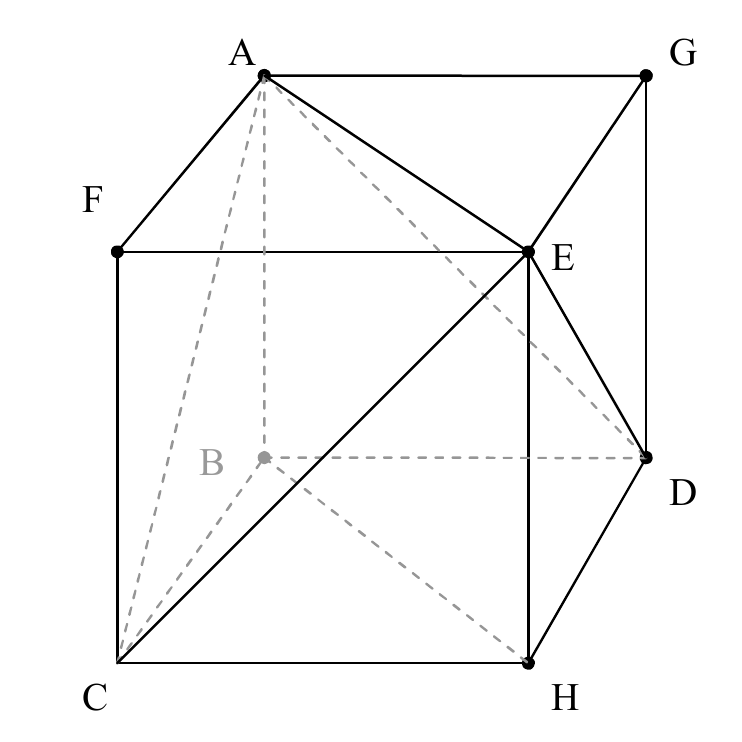}
\includegraphics[height=7.5cm]{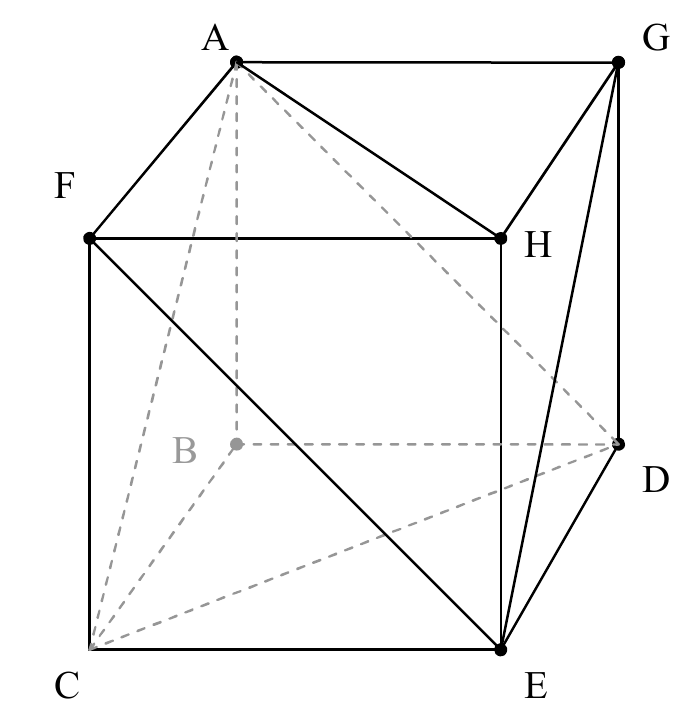}}\caption{Graphs $G_{8,3}$ and $G_{8,4}$} \label{p20} \end{figure}

\begin{figure}[ht]  \center{\includegraphics[height=7.5cm]{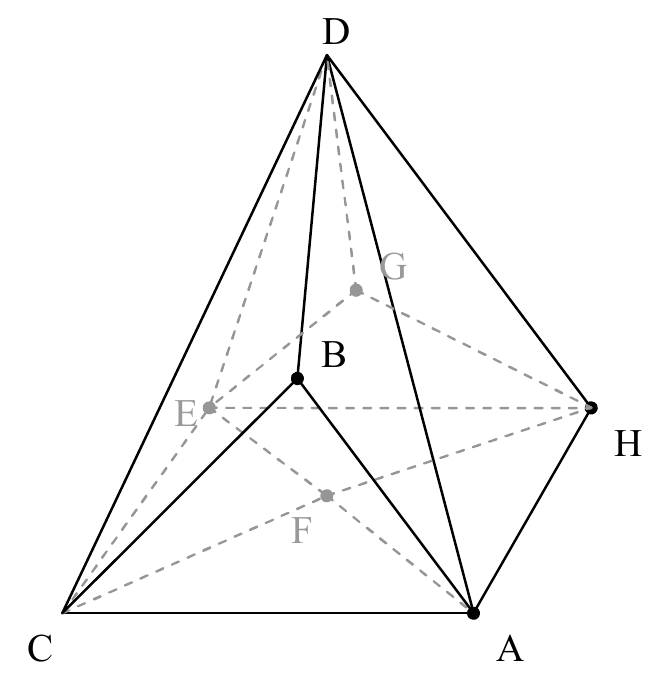}
\includegraphics[height=7.5cm]{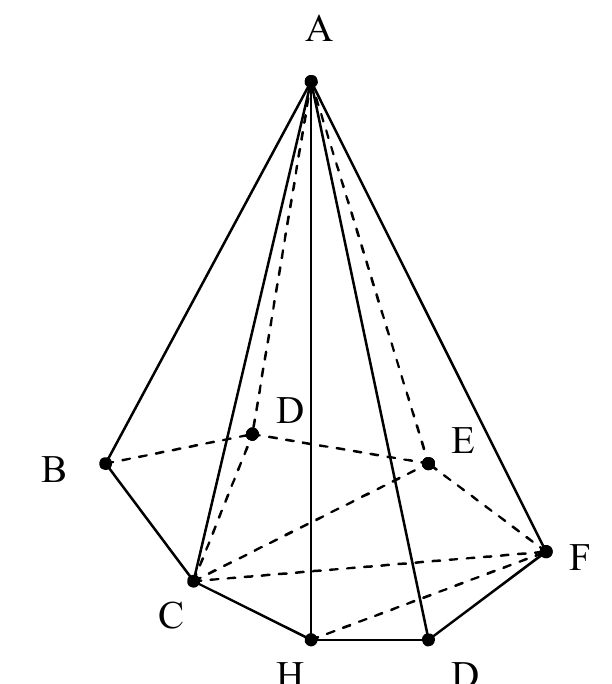}}\caption{Graphs $G_{8,5}$ and $G_{8,6}$} \label{p23} \end{figure}

\hfill
\newpage
\begin{figure}[ht]  \center{\includegraphics[height=7.5cm]{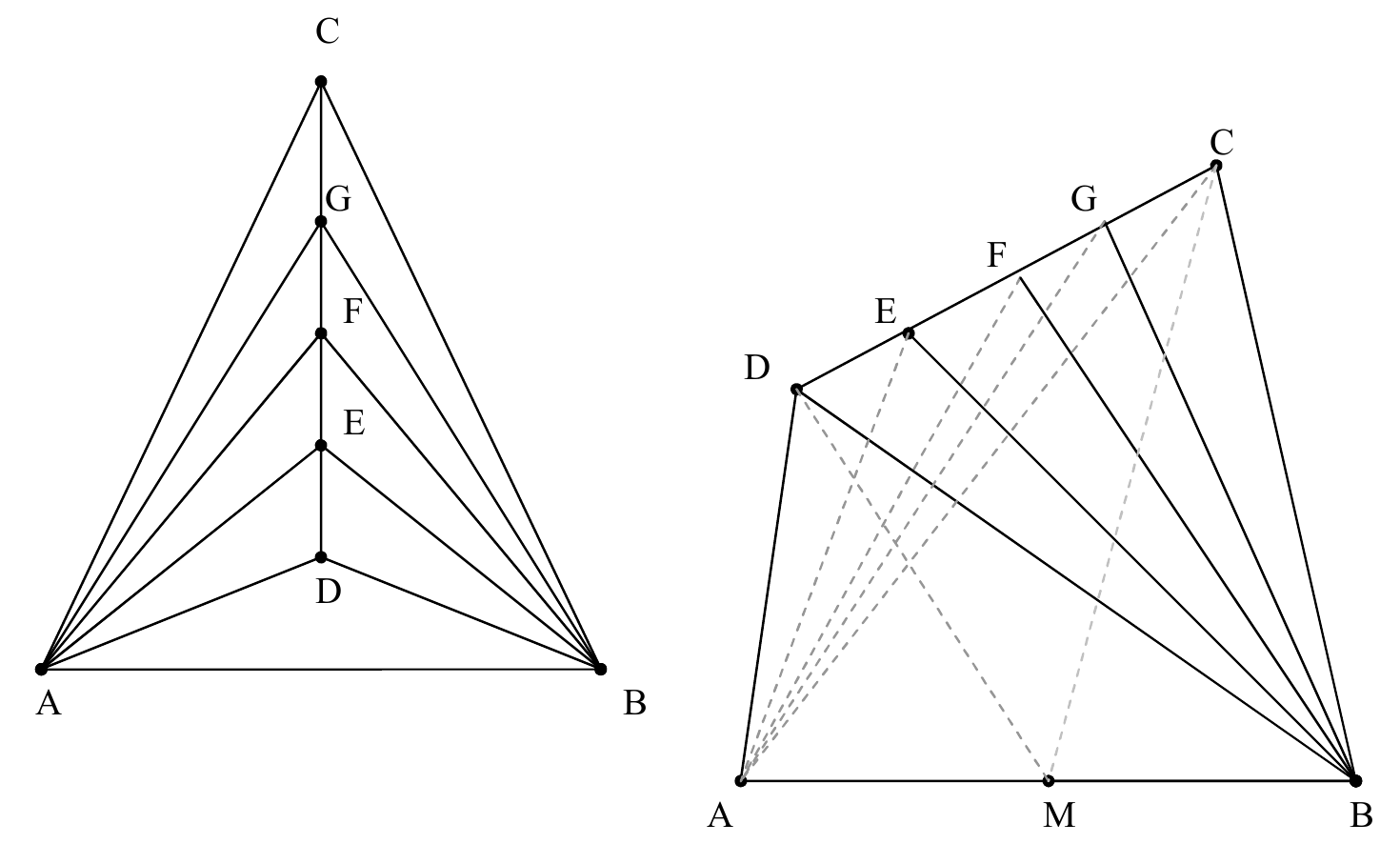}}\caption{Graph $G_{8,8}$} \label{p22} \end{figure}


\begin{figure}[ht]  \center{\includegraphics[height=6.5cm]{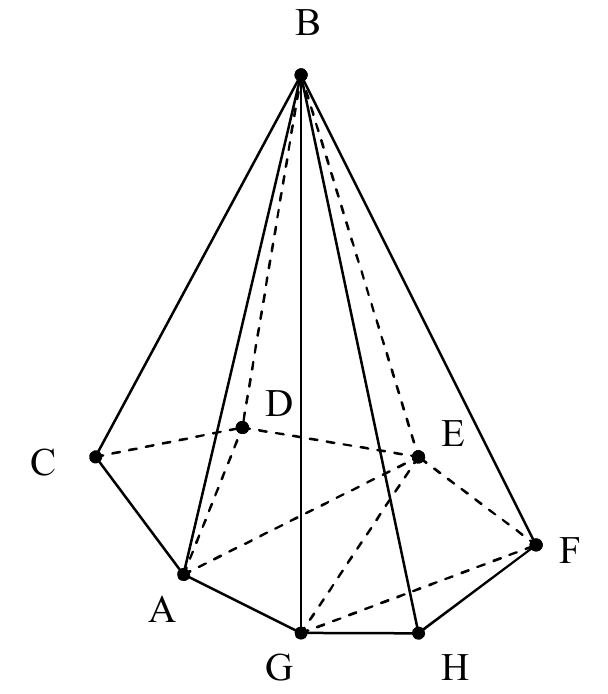}}\caption{Graph $G_{8,9}$} \label{p24} \end{figure}
\begin{figure}[ht]  \center{\includegraphics[height=7.5cm]{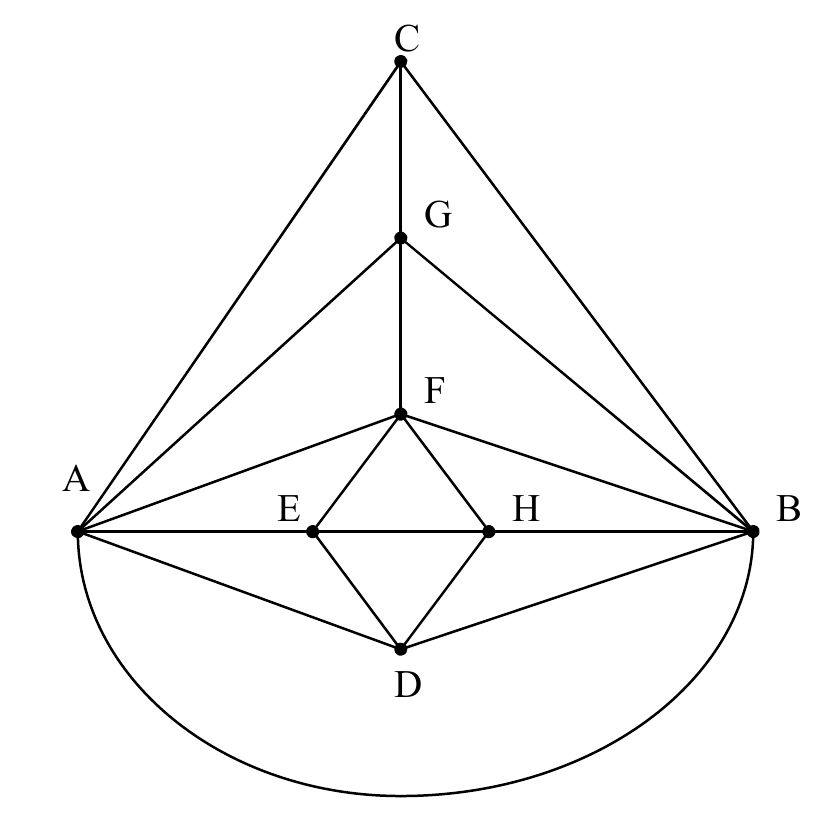}}\caption{Graph $G_{7,4}$ and the plane of symmetry.} \label{p25} \end{figure}



\ \ Consider the graph G \textsubscript{7.4 }(Fig. \ref{p25}). As we see, there are two planes of
symmetry: { S} DM ( M is the middle of AB ) and ABF. Then
realized that in this time it is enough to look at the next and only the next attack:

(10) BD ${\approx}$ BDE;

(11) BE;

(12) B.F.;

(13) AB;

(14) DE ${\approx}$ EF ${\approx}$ ABD;

(15) BEF;

(16) ABDEFA;

(17) ADEFGA ${\approx}$ AEFGCA.

Non-isomorphic in this time are 5 splits: G \textsubscript{8.10 }, G \textsubscript{8.11 }, G \textsubscript{8.13 }, G
\textsubscript{8.14 }, G \textsubscript{8.16 }(div. Fig. 2 6 - 30 ). For other cuts, maybe G \textsubscript{8.10
}${\approx}$ G \textsubscript{8.8 }, G \textsubscript{8.15 }${\approx}$ G \textsubscript{8.1 }, G \textsubscript{8.17
}${\approx}$ G \textsubscript{8.6 }.
\newpage

\begin{figure}[ht]  \center{\includegraphics[height=7.5cm]{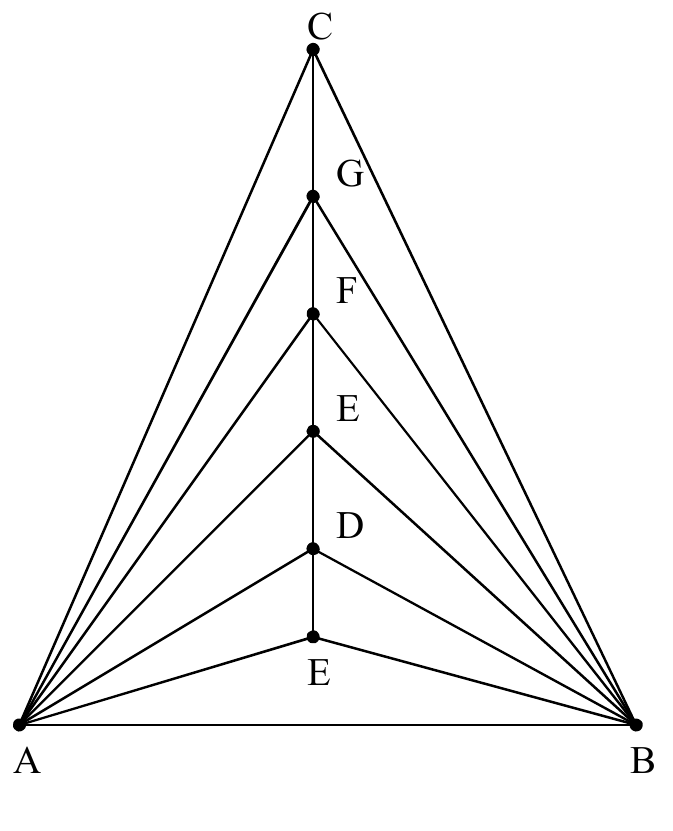}
\includegraphics[height=7.5cm]{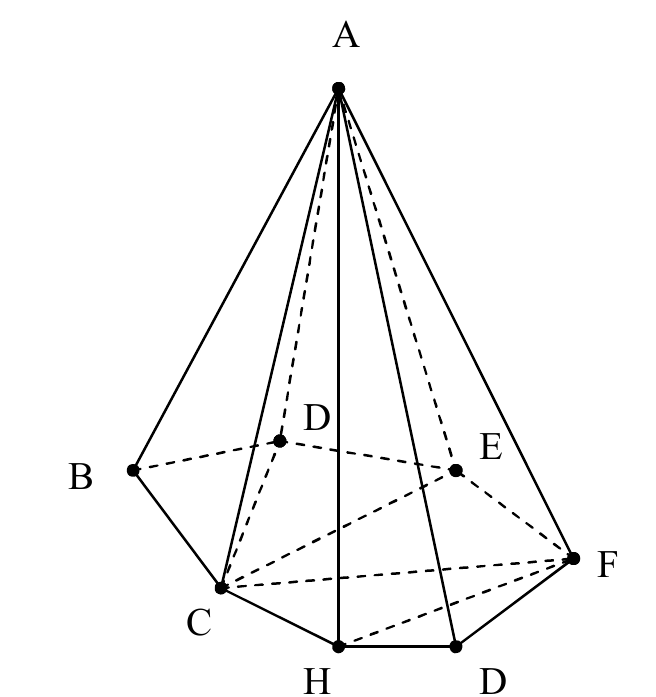}}\caption{Graphs $G_{8,10}$ and $G_{8,16}$} \label{p27} \end{figure}
\hfill
\newpage

\begin{figure}[ht]  \center{\includegraphics[height=7.5cm]{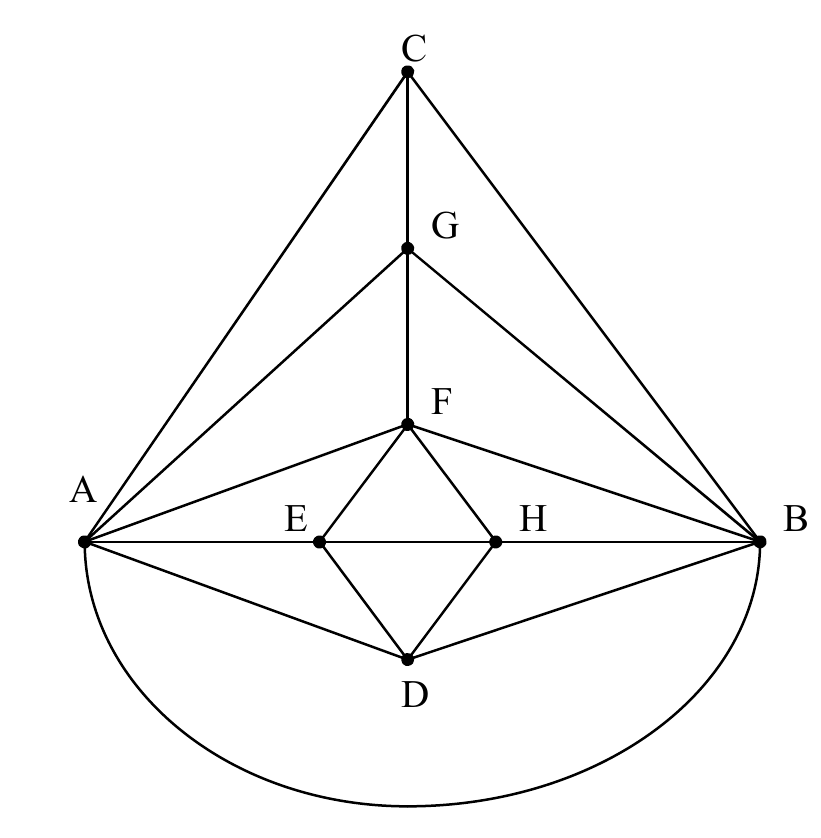}
\includegraphics[height=7.5cm]{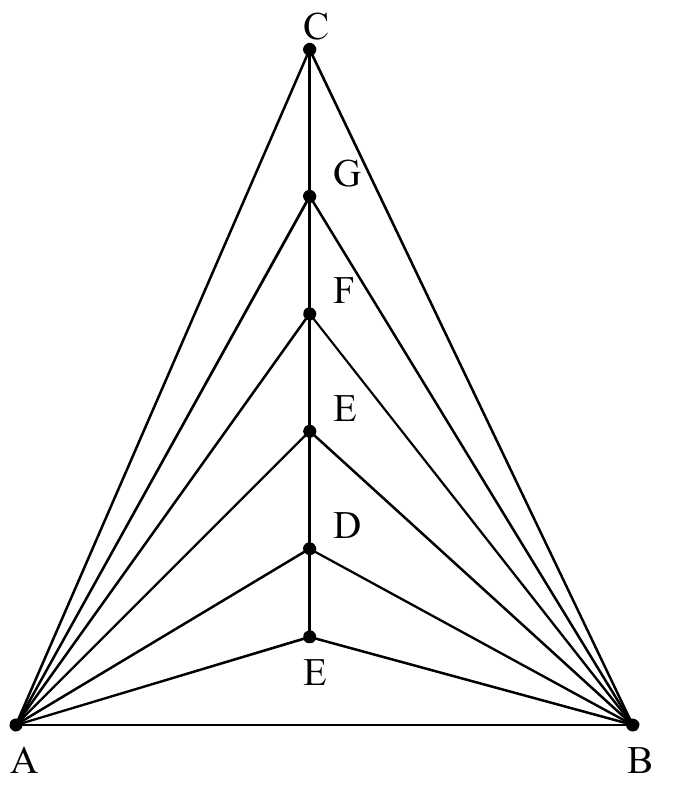}}
\caption{Graphs $G_{8,11}$ and $G_{8,14}$} 
\label{p29} 
\end{figure}


\begin{figure}[ht]  \center{\includegraphics[height=7.5cm]{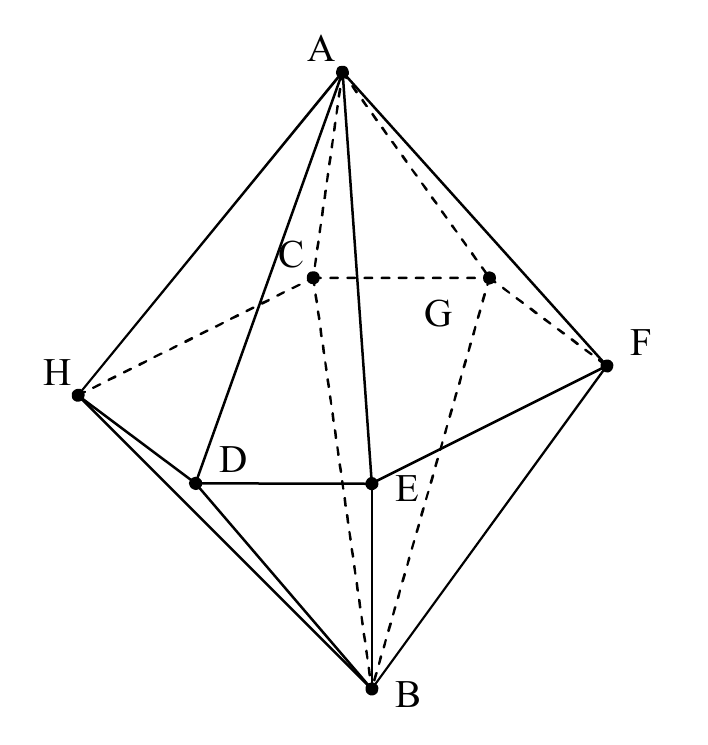}}\caption{Graph $G_{8,13}$} \label{p30} \end{figure}


\bigskip

Let's consider graph (rozbittya) G \textsubscript{7.3 }(Fig. 15, 30). As we see, there are two planes
of symmetry: AMN ( de M is the middle of the DE , N is the middle of the BF ) and ACG . The wanderers from the circles
realized that in this time it is enough to look at the onset and only the onset of the edges:

(1 8 ) BE;

( 20 ) EF;

( 2 2) AG;

( 24 ) BED;

( 26 ) AEG;

( 28 ) DEGACD;

( 30 ) EDBCAE; 

(1 9 ) DE;

( 2 1) EA;

(2 3 ) GE;

( 25 ) ADE;

( 27 ) EGF;

( 29 ) EABCDE;

( 3 1) EBCADE.

Only 2 divisions are not isomorphic in this case: G \textsubscript{8.27 }and G \textsubscript{8.28} (Div. Fig. 31-32 ) .
For other cuts, maybe G \textsubscript{8.1 8 }${\approx}$ G \textsubscript{8, 9 }, G \textsubscript{8.1 9 }${\approx}$
G \textsubscript{8.26 }${\approx}$ G \textsubscript{8.1 0 }, G \textsubscript{8, 20 }${\approx}$ G \textsubscript{8, 22
}${\approx}$ G \textsubscript{8, 11 }, G \textsubscript{8, 21 }${\approx}$ G \textsubscript{8.6 }, G
\textsubscript{8.23 }${\approx}$ G \textsubscript{8.16 }, G \textsubscript{8.24 }${\approx}$ G \textsubscript{8.3 0
}${\approx}$ G \textsubscript{8.2 }, G \textsubscript{8.25 }${\approx}$ G \textsubscript{8.3 0 }${\approx}$ G
\textsubscript{8.1 }, G \textsubscript{\_ 8.29 }${\approx}$ G8.27 . \textsubscript{\_}


\begin{figure}[ht]  \center{\includegraphics[height=6.5cm]{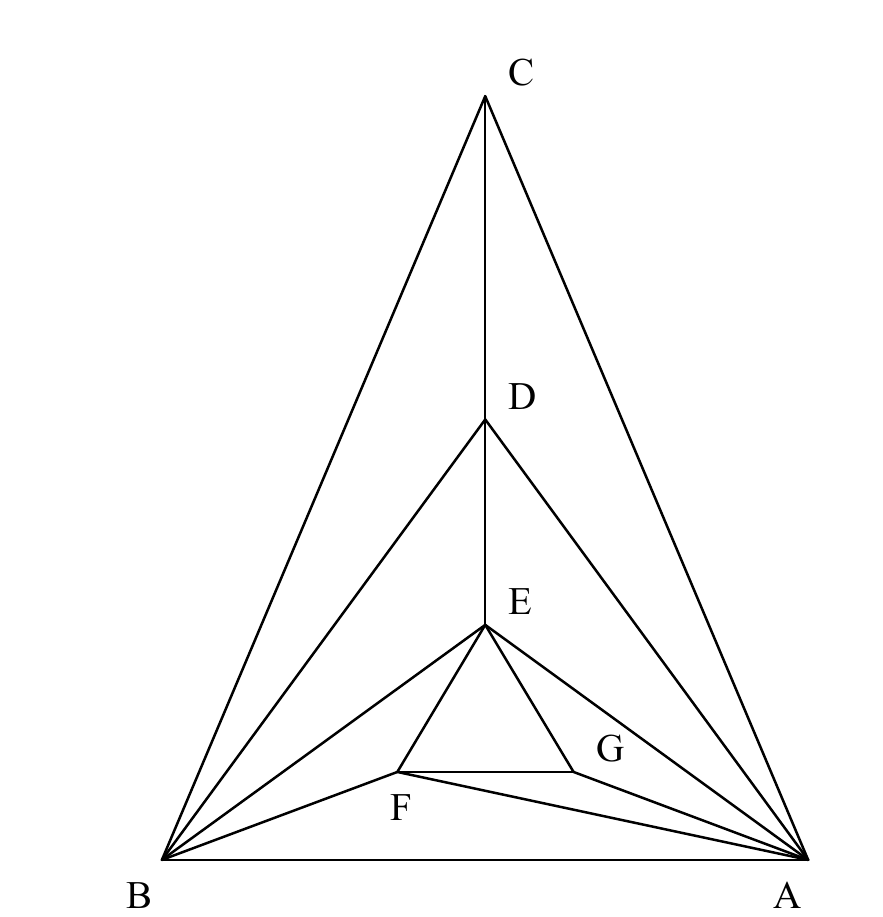} \includegraphics[height=6.5cm]{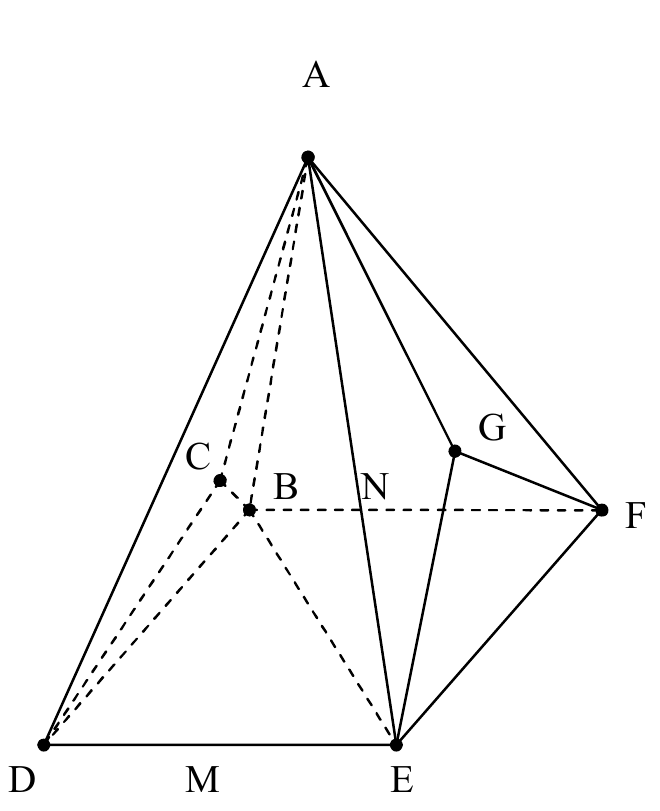}}\caption{Graph $G_{7,3}$ Graph G \textsubscript{7.3 }( p. fig. 15) and the plane of symmetry.} 

\label{p32} \end{figure}



\begin{figure}[ht]  \center{\includegraphics[height=6.5cm]{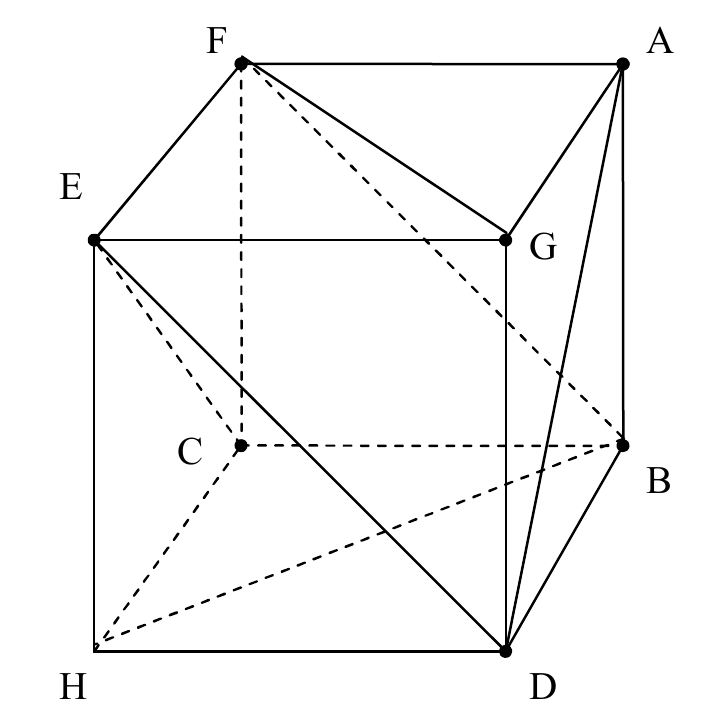}
\includegraphics[height=6.5cm]{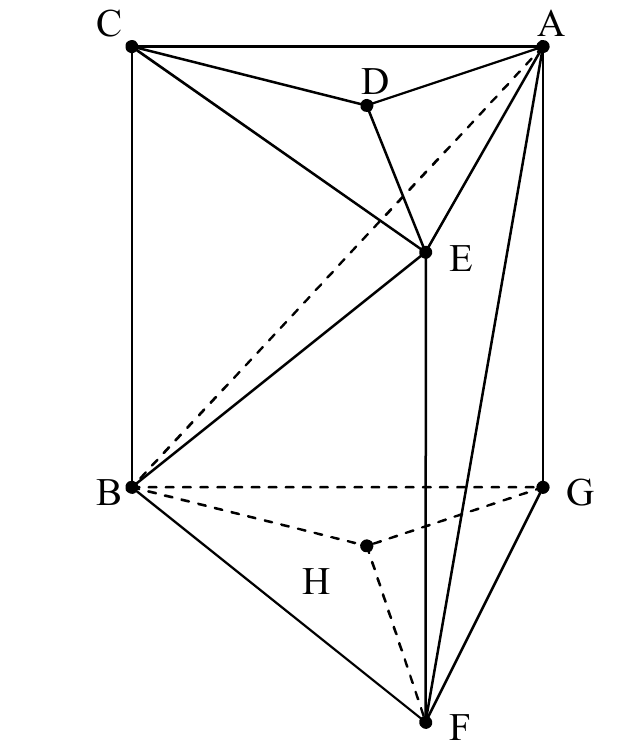}}\caption{Graphs $G_{8,28}$ and $G_{8,27}$} \label{p34} \end{figure}


Consider the  graph  G \textsubscript{7.1 }.
We see, we see all central symmetry do. Walkers from the mirroring of symmetry realized that in this time it is enough to look at the offensive and less than
the offensive edges:

(32) AB;

(34) FG;

(36) ABD;

(38) AGB;

(40) BECFGB;

(42) BCEGAB; 

(33) AG;

(35) AD;

(37) AGF;

(39) FEG;

(41) BEGADB;

(43) BGFADB.

Non-isomorphic in this case, there is only one difference: G \textsubscript{8.39} (div. Fig. 34). For other cuts, maybe
G \textsubscript{8.32 }${\approx}$ G \textsubscript{8.41 }${\approx}$ G \textsubscript{8.43 }${\approx}$ G
\textsubscript{8.28 }, G \textsubscript{8.33 }${\approx}$ G \textsubscript{8.34 }${\approx}$ G \textsubscript{8.6 }, G
\textsubscript{8.35 }${\approx}$ G \textsubscript{8.36 }${\approx}$ G \textsubscript{8.42 }${\approx}$ G
\textsubscript{8.11 }, G \textsubscript{8.37 }${\approx}$ G \textsubscript{8.9 }, G \textsubscript{8.38 }${\approx}$ G
\textsubscript{8.4 }, G \textsubscript{8.40 }${\approx}$ G \textsubscript{8.4 }.

\begin{figure}[ht]  \center{\includegraphics[height=6.5cm]{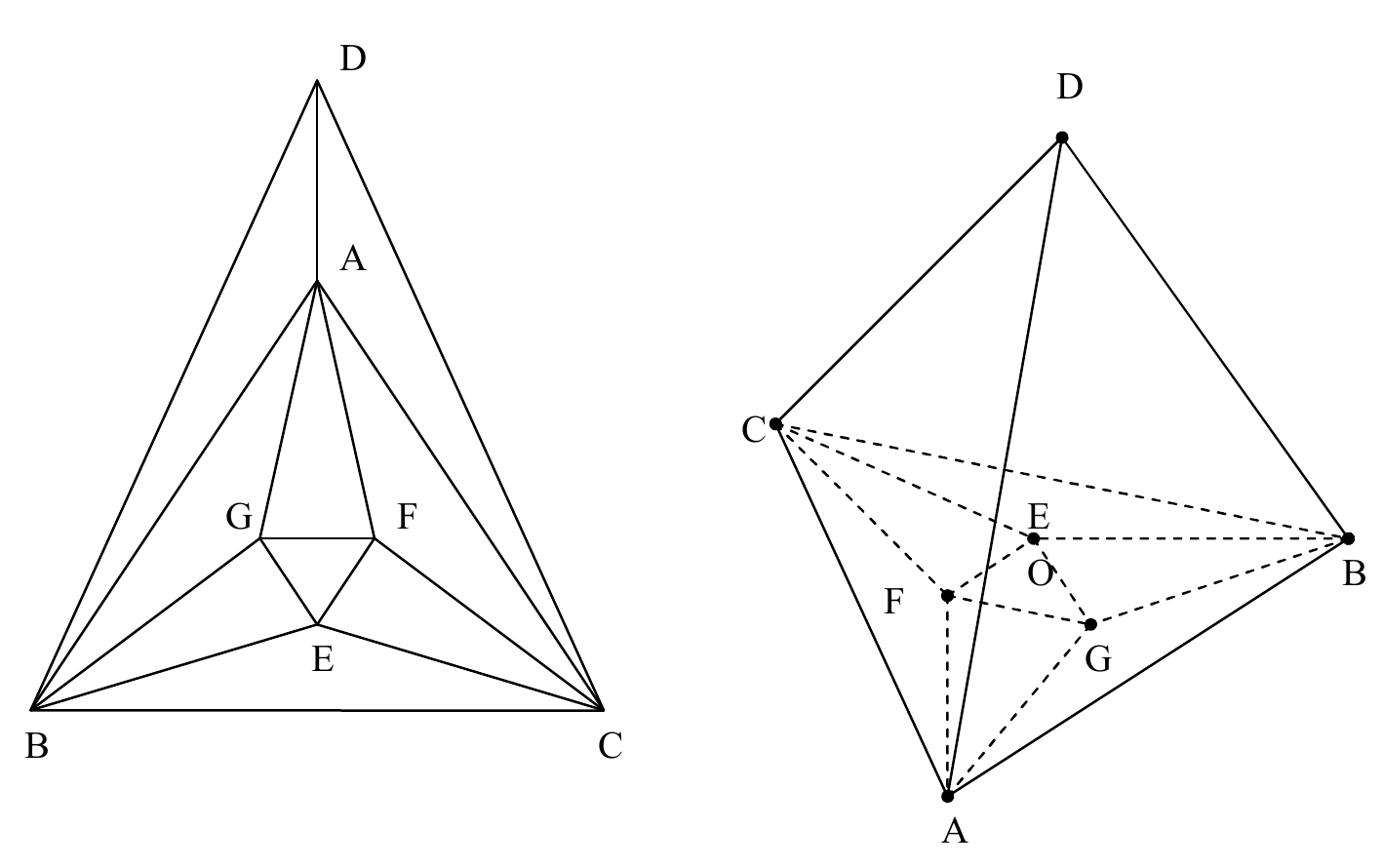}}\caption{Graph $G_{7,1}$ and all symmetry.} \label{p35} \end{figure}



\begin{figure}[ht]  \center{\includegraphics[height=6.5cm]{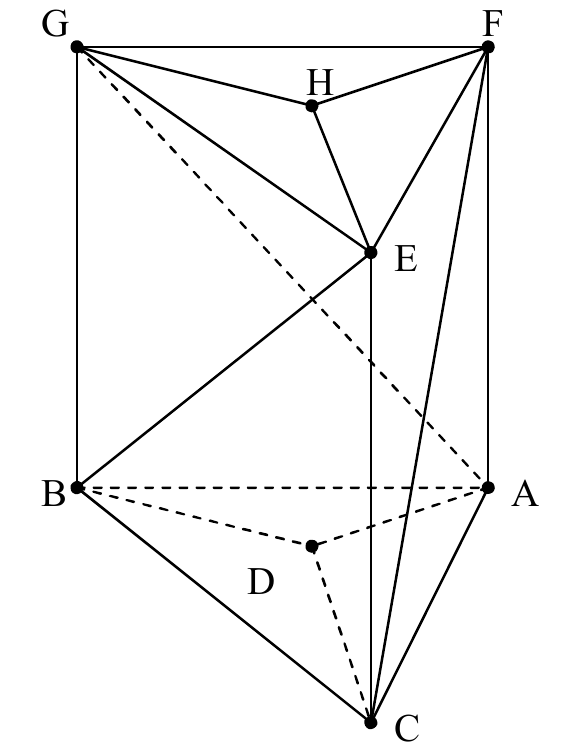}}\caption{Gpaph $G_{8,39}$} \label{p36} \end{figure}

\begin{figure}[ht]  \center{\includegraphics[height=6.5cm]{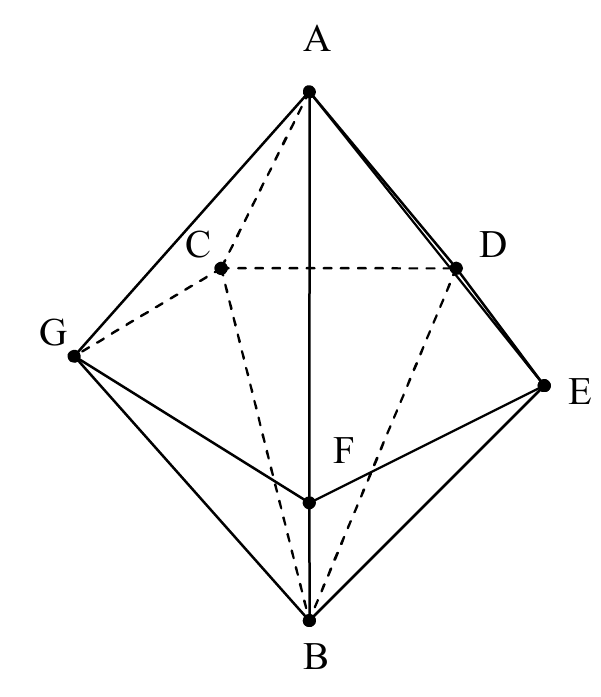}}\caption{Gpaph $G_{7,2}$} \label{p37} \end{figure}



\bigskip

Consider the graph G \textsubscript{7.2 }(Fig. 14, 35). We see, it is centrally symmetrical.
Walkers from the mirroring of symmetry realized that in this time it is enough to look at the offensive and less than
the offensive battle:

( 43) AC;

( 4 4) CD;

( 45 ) ACD;

( 46 ) BCDEFG;

(4 7 ) ADBFGA.

\ \ Mustaches of strife appeared to be isomorphic to look at; and itself: G \textsubscript{8.43 }${\approx}$ G
\textsubscript{8.46 }${\approx}$ G \textsubscript{8.28 }, G \textsubscript{8.44 }${\approx}$ G \textsubscript{8.13 }, G
\textsubscript{8.45 }${\approx}$ G \textsubscript{8.47 }${\approx}$ G \textsubscript{8.6 }.

\ \ In this rank, the task was set to rozv ' yazana povnistyu; $\mu $(8)=14. In Table \ref{t2}, there is a list of obsessions
of battles G \textsubscript{8 }from the grades of their vertices. Sliding respect for those graphs G \textsubscript{8.4
}and G \textsubscript{8.27 }, appeared to be non-isomorphic, not wondering at the same set of steps of vertices.

\begin{table}[ht]
	\centering
		\begin{tabular} {|c|c|c|c|c|c|c|c|c|}
		
\hline
G\textsubscript{8,3} &
3 &
3 &
3 &
3 &
6 &
6 &
6 &
6\\ 
\hline
G\textsubscript{8,2} &
3 &
3 &
3 &
4 &
5 &
6 &
6 &
6\\
\hline
G\textsubscript{8,1} &
3 &
3 &
3 &
4 &
5 &
5 &
6 &
7\\
\hline
G\textsubscript{8,14} &
3 &
3 &
4 &
4 &
4 &
4 &
7 &
7\\ 
\hline
G\textsubscript{8,10} &
3 &
3 &
4 &
4 &
4 &
5 &
6 &
7\\
\hline
G\textsubscript{8,16} &
3 &
3 &
4 &
4 &
5 &
5 &
5 &
7\\
\hline
G\textsubscript{8,4} &
3 &
3 &
4 &
4 &
5 &
5 &
6 &
6\\
\hline
G\textsubscript{8,27} &
3 &
3 &
4 &
4 &
5 &
5 &
6 &
6\\
\hline
G\textsubscript{8,9} &
3 &
3 &
4 &
5 &
5 &
5 &
5 &
6\\
\hline
G\textsubscript{8,39} &
3 &
3 &
5 &
5 &
5 &
5 &
5 &
5\\
\hline
G\textsubscript{8,11} &
3 &
4 &
4 &
4 &
4 &
5 &
6 &
6\\
\hline
G\textsubscript{8,6} &
3 &
4 &
4 &
4 &
5 &
5 &
5 &
6\\
\hline
G\textsubscript{8,13} &
4 &
4 &
4 &
4 &
4 &
4 &
6 &
6\\
\hline
G\textsubscript{8,28} &
4 &
4 &
4 &
4 &
5 &
5 &
5 &
5\\
\hline

\end{tabular}

\caption{Degree of vertexes} 
\label{t2}
\end{table}

As we see in Table 1, the sum of the steps of all vertices for any kind of breakup is good ( n {}-2) ${\times}$6,
which, no matter how you show, is true for any breakout of spheres from help n dot. In Table 2, obviously, there are
far from all combinations of eight natural numbers from 3 to 7, the sum of which is more than 36. In the development of
a simple algorithm for re-verification of such a possibility, the problem of searching for solving the sphere, for
whatever n , could be taken into account. It seems that this problem will be examined in our further investigations.

Realization of the method we used for small n is possible for $n <12$. The solution of the
problem of splitting the sphere into triangles with additional more than 8 vertices, obviously, is due to the solving.

\bigskip

\section{Tricolor coloring of 1-skeleton}

The three-color coloring of a 1-skeleton is the coloring of its edges in three colors, in which the boundary of each face contains edges of all three colors. Note that if there is a vertex of degree 3, then all edges incident to it are painted in different colors. Two colorings will be considered the same (equivalent) if one can be obtained from the other by renaming the colors.

\begin{figure}[ht]  \center{\includegraphics[height=5.5cm]{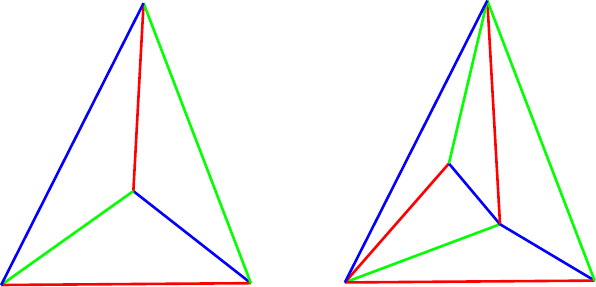}}\caption{Unique coloring of $G_4$ and $G_5$} \label{p41} \end{figure}


Let us show that the graph $G_4$ has a unique coloring. The three outer edges can be colored in a single way, with the accuracy of the substitution of colors. Let's fix this coloring. Then, for each external vertex, the color of the internal edge incident to it is determined by the colors of two adjacent external edges. So we get the only possible coloring of the 1-skeleton.


Let us show the uniformity of the coloring of graph $G_5$. We use the following colors: red, green, blue. Let's fix the colors of the sides of the outer triangle again. Let the lower edge be red and the right edge be green, and let the lower right vertex $A$ have degree three. Then the third edge coming out of it is blue. Then the coloring of the two inner triangles with vertex $A$ is unambiguously given. Three edges coming from one vertex remain uncolored. The process of their coloring is the same as for the $G_4$.

\begin{figure}[ht] 
\center{\includegraphics[height=7.0cm]{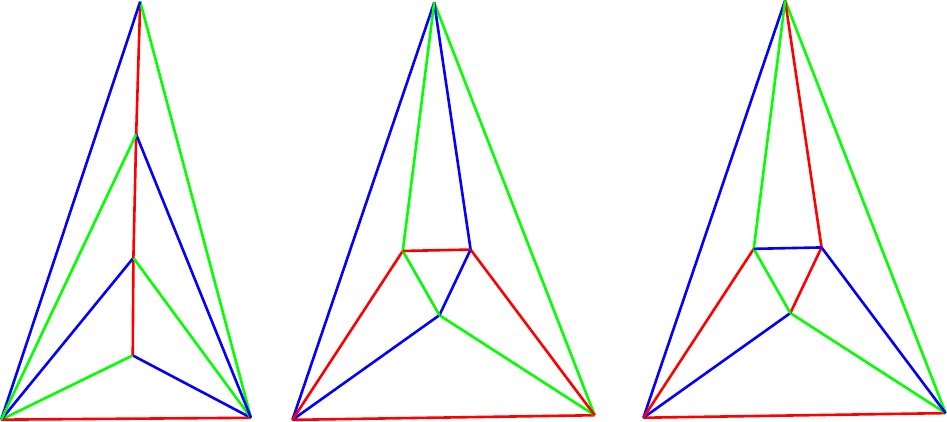}}
\caption{Unique colorings  of $G_{6,1}$ and   two colorings of $G_{6,2}$} \label{p42} \end{figure}


By analogy with $G_5$, starting from the upper vertex, we get a unique coloring of $G_{6,1}$. For $G_{6,2}$, there are two colorings shown in fig.\ref{p42} center and left. These colorings are not equivalent, because in the first coloring the edges of each color form a cycle of length 4, and in the second coloring the union of red (or similarly blue) edges is a disjoint set.

\begin{figure}[ht]  
\center{\includegraphics[height=6.5cm]{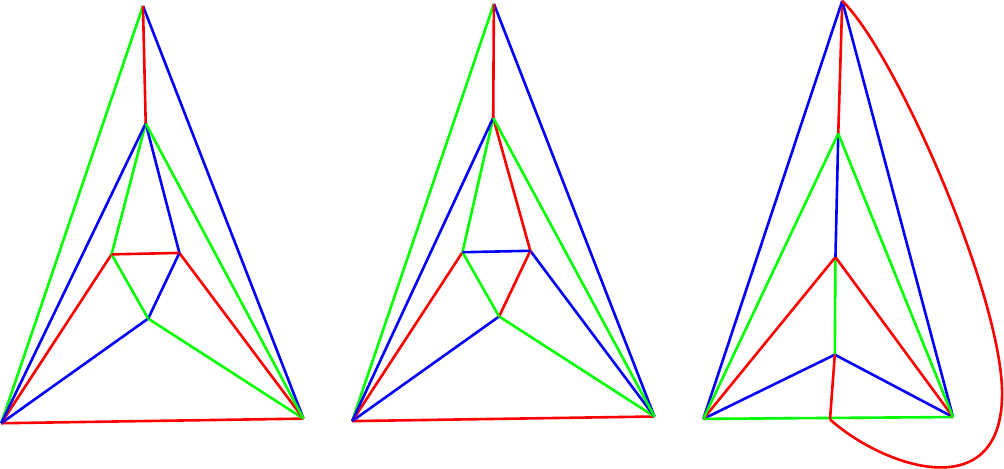}}
\caption{Two colorings  of $G_{7,1}$ and   unique coloring of $G_{7,2}$} \label{p43} \end{figure}


Two colorings of $G_{6,2}$ give rise to two colorings of $G_{7,1}$. They are not equivalent because in the first graph the edges of each color have two components of 1 and 4 edges, and in the second graph the components have 2 and 3 edges. For other graphs with 7 vertices, the coloring is unique and shown in the remaining figures.

\begin{figure}[ht]  \center{\includegraphics[height=7.0cm]{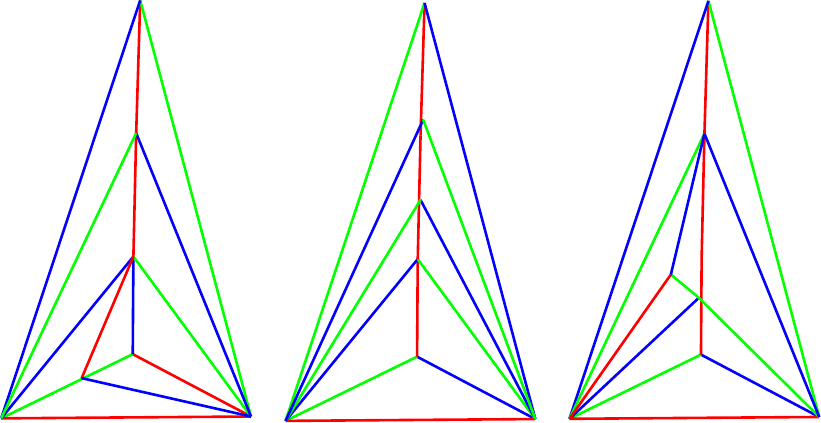}}
\caption{Unique colorings  of $G_{7,3}$,$G_{7,4}$,$G_{7,5}$ } \label{p44} \end{figure}


If we summarize all of the above, we get the following 
\begin{theorem} For one triangulation of a sphere with 6 vertices and one triangulation with 7 vertices, there are two colorings, for the rest of the triangulations with no more than 7 vertices, there is a single coloring.
\end{theorem}



\newpage

\bibliographystyle{plain}
\bibliography{abp}

\end{document}